\documentclass[11pt]{article}

\usepackage[T1,T2A]{fontenc}
\usepackage[utf8]{inputenc}
\usepackage{amsfonts}
\usepackage{amssymb}
\usepackage{amsmath}
\usepackage{epsfig}
\usepackage{color}
\usepackage{theorem}
\usepackage{tikz}

\newtheorem{theorem}{Theorem}[section]
\newtheorem{lemma}[theorem]{Lemma}
\newtheorem{corollary}[theorem]{Corollary}
\newtheorem{proposition}[theorem]{Proposition}

{\theorembodyfont{\rmfamily} 

}

\usepackage[color links=false, backref=page]{hyperref} 

\def\Q{\mathbb{Q}}

\def\F{\mathbb{F}}
\def\R{\mathcal{R}}
\def\O{\mathcal{O}}
\def\H{\mathbb{H}}
\def\K{\mathcal{K}}
\def\Fr{\mathcal{F}}
\def\M{\mathcal{M}}
\def\N{\mathcal{N}}
\def\S{\mathcal{S}}
\def\SL{\rm{SL}}
\def\PSL{\rm{PSL}}
\def\A{\rm{Aut}}

\title{Hole operations on Hurwitz maps}
\date{\empty}
\author{G\'abor G\'evay and Gareth A. Jones}

\begin{document}

\medskip

\maketitle

\centerline{Dedicated to Marston Conder in celebration of his 65th birthday.}

\bigskip

\begin{abstract}
For a given group $G$ the orientably regular maps with orientation-preserving automorphism group $G$ are used as the vertices of a graph 
$\O(G)$, with undirected and directed edges showing the effect of duality and hole operations on these maps. Some examples of these graphs 
are given, including several for small Hurwitz groups. For some $G$, such as the affine groups ${\rm AGL}_1(2^e)$, the graph $\O(G)$ is 
connected, whereas for some other infinite families, such as the alternating and symmetric groups, the number of connected components is 
unbounded.
\end{abstract}

\noindent{\bf MSC classification:} 05C10 (primary), 
20B25 (secondary). 
 
\medskip

 \noindent{\bf Key words and phrases:} regular map, orientably regular, duality, hole operation, Hurwitz group.

\section{Introduction}\label{sec:intro}

For any group $G$ let $\R(G)$ and $\O(G)$ denote the sets of (isomorphism classes of) regular or orientably regular maps $\M$ with automorphism 
group ${\rm Aut}\,\M$ or orientation-preserving automorphism group ${\rm Aut}^+\M$ isomorphic to $G$. For a given finite group $G$ these two 
sets are both finite since, in each case, their elements correspond bijectively to the orbits of ${\rm Aut}\, G$ on certain generating triples for $G$. 
These maps are related to each other by the application of certain operations, described by Wilson in~\cite{Wil79}, such as duality $D$, Petrie duality 
$P$ (in the case of regular maps) and the hole operations $H_j$ for $j$ coprime to the valency of the vertices, since all of these preserve the group 
$G$. One way of understanding how the various maps associated with $G$ are related to each other is to regard these sets $\R(G)$ and $\O(G)$ as 
the vertex sets of graphs, with directed or undirected edges indicating the actions of these operations. For example, one can consider whether or not 
these graphs are connected, and if not, what invariants might be used to distinguish maps in different components. 
We will describe the maps $\M$ and graphs $\O(G)$ arising for some particularly interesting groups, such as a few of the 
smallest Hurwitz groups. We will show that for the symmetric and alternating groups, the number of connected components of $\O(G)$ is unbounded,
building on work of Marston Conder~\cite{Con80} on alternating Hurwitz groups in the latter case.


\section{Orientably regular maps}\label{sec:orreg}

If $\M$ is an orientably regular map then its arcs (directed edges) can be identified with the elements of a group $G$, so that its orientation-preserving monodromy 
and automorphism groups are identified with the right and left regular representations of $G$.

Let $\M$ have type $\{p,q\}$, meaning that its faces are $p$-gons and its vertices have valency $q$. (Here we follow Coxeter's notation, see~\cite{CM}, for example.) 
Then $G$ has generators $x,y, z$ satisfying
\[x^q=y^2=z^p=xyz=1,\]
where $x$ and $y$ act (as monodromy permutations) by rotating all arcs around their incident vertices, following the orientation, and by reversing them, so that $z$ 
rotates them around faces. Conversely any such generating triple for $G$ defines a map $\M$ of this type, with isomorphic maps corresponding to triples equivalent 
under ${\rm Aut}\,G$. Of course the triple, and hence the map, is uniquely determined by any two of its members, so we will often restrict attention to the generating 
pair $x, y$.

Map duality $D$, transposing vertices and faces, replaces the triple $(x,y,z)$ with $(z, y, (zy)^{-1}=x^y)$, corresponding to the dual map $D(\M)$. 
Then $D^2$ sends $(x,y,z)$ to $(x^y,y,z^y)=(x,y,z)^y$, giving a map isomorphic to $\M$.

\begin{figure}[!h]
\begin{center}
\includegraphics[width=0.85\textwidth] {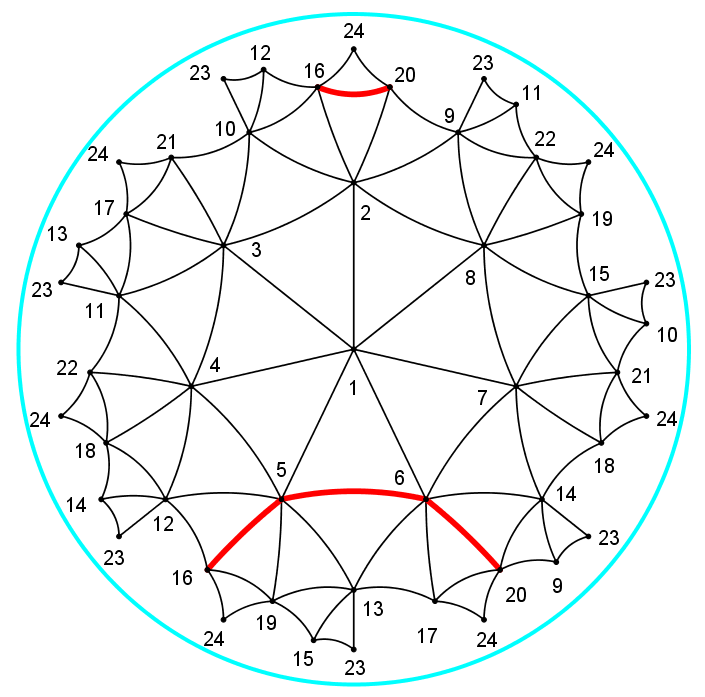}
\caption{Klein's map $\K$ of type $\{3,7\}$ and genus $3$ represented in Poincar\'e's disc model of the hyperbolic plane.
A 3-hole $(16, 5, 6, 20,16)$ is highlighted.}
\label{fig:Klein}
\end{center}
\end{figure}

If $j$ is coprime to $q$ then another orientably regular map $\M_j=H_j(\M)$, of type $\{p',q\}$ for some $p'$ and with the same orientation-preserving 
automorphism group $G$ as $\M$, can be found by applying the hole operation $H_j$, described by Wilson in~\cite{Wil79}, to $\M$. This map embeds 
the same graph as $\M$, but the rotation $x$ of arcs around each vertex is replaced with $x^j$, so that the faces, and hence the underlying surface, 
may be changed. In terms of group theory, $H_j(\M)$ corresponds to the generating triple $(x^j, y, z_j)$ for $G$, where $z_j:=(x^jy)^{-1}$. Thus 
the boundary of each face of $H_j(\M)$, given the orientation of the underlying surface,
follows a $j$-{\em hole\/} of $\M$, at each vertex $v$ taking the $j$-th incident edge, following the reverse 
orientation around $v$, rather than the first; its length is the order of the element $z_j$.
(Note that a $j$-hole, taken in the reverse direction, becomes a $(q-j)$-hole.)

\medskip

\noindent{\bf Example} Figure~\ref{fig:Klein} shows a $3$-hole $( 16, 5, 6, 20, 16)$ of length $4$ in
Klein's map $\K$ of type $\{3,7\}$ and genus $3$ with $G={\rm PSL}_2(7)$ (see Section~\ref{sec:PSL_2(7)} for details).

\medskip

Clearly $H_j=H_k$ whenever $j\equiv k$ mod~$(q)$; since $H_j\circ H_k=H_{jk}$ for all $j$ and $k$ coprime to $q$, these operations give a 
representation of the group $U_q$ of units mod~$(q)$. The operation $H_{-1}$ sends each map to its mirror-image, so an orientably regular 
map $\M$ is regular (has a flag-transitive full automorphism group ${\rm Aut}\,\M$) if and only if $H_{-1}(\M)\cong \M$; thus the action on 
regular maps gives a representation of the group $U_q/\{\pm 1\}$.    \newpage

Another useful map operation is Petrie duality. A {\em Petrie polygon\/} in a map $\M$ is a closed zig-zag path turning alternately first right and 
first left at successive vertices. The {\em Petrie dual\/} $P(\M)$ of $\M$ embeds the same graph as $\M$, but the faces of $\M$ are replaced 
with new faces bounded by the Petrie polygons. If $\M$ is orientably regular their length $r$ (the {\em Petrie length\/} of $\M$) is twice the 
order of the commutator $[x,y]=x^{-1}y^{-1}xy=x^{-1}yxy$, so that $P(\M)$ has type $\{r,q\}$. As in~\cite{Con09, CM}, we will often refer to
the {\em extended type\/} $\{p,q\}_r$ of $\M$. For instance, in Figure~\ref{fig:Klein} the $3$-holes $(16, 5, 6, 20,16)$ and $(24, 17, 13, 19, 24)$ 
in $\K$ enclose a Petrie polygon of length $r=8$. (More generally, whenever $p=3$ each Petrie polygon is enclosed by two $3$-holes of length $r/2$ 
in this way.)

Given any group $G$, one can regard $\O(G)$ as the vertex set of a graph
(also denoted by $\O(G)$) by adding edges showing the actions 
of the operations $D$ and $H_j$ (for selected $j$); they will be undirected for cycles of length $2$, e.g.~for involutions such as $D$, and 
directed in the case of longer cycles for hole operations $H_j$. Loops, corresponding to invariant maps, will be omitted. In particular, we 
will use dashed and dotted undirected edges for $D$ and $H_{-1}$, and unbroken edges for other hole operations. (We will not represent 
the action of Petrie duality $P$ in a similar way, since although $P$ preserves automorphism groups, it does not always preserve orientability: 
for example, if $\M$ is the tetrahedral map on the sphere, then $P(\M)$ is the antipodal quotient of the cube, on the real projective plane; 
similarly, $P(\K)$ is a non-orientable regular map of genus $41$.)

Provided its genus $g$ is neither too large nor less than $2$, each map $\M\in\O(G)$ can be located in Conder's computer-generated lists of 
regular or chiral orientable maps in~\cite{Con09}. Each entry there has the form R$g.n$ or C$g.n$ respectively, where $g$ denotes the genus, 
and $n$ denotes the $n$th entry in the list of maps of that genus, ordered lexicographically by their type $\{p,q\}$, with $p\le q$. Within each 
list, each entry refers either to a dual pair of maps of types $\{p,q\}$ and $\{q,p\}$ or to a single self-dual map with $p=q$. If $p<q$ we will 
denote the maps of type $\{p,q\}$ and $\{q,p\}$ by X$g.n$a and X$g.n$b respectively, where X is R or C; if $p=q$ and the entry denotes a 
non-isomorphic dual pair we will assign the labels X$g.n$a and X$g.n$b arbitrarily, whereas a single self-dual map will be denoted simply by 
X$g.n$. Similar conventions will apply to the list of non-orientable regular maps, denoted by N$g.n$ in~\cite{Con09}.


\subsection{Some simple examples of graphs $\O(G)$}\label{sec:small}

If $G$ is the alternating group ${\rm A}_4$ then $\O(G)$ contains only the tetrahedral map $\{3,3\}$ of genus $0$.

If $G$ is a dihedral group ${\rm D}_n$ for some $n>2$ then $\O(G)$ consists of the dual pair of maps $\{2,n\}$ and $\{n,2\}$ of genus $0$, 
whereas if $n=2$, so that $G={\rm D}_2\cong {\rm V}_4$, then $\O(G)$ contains only the self-dual map $\{2,2\}$ of genus~$0$. In either 
case the operations $H_j$ for $j$ coprime to $n$ act trivially on $\O(G)$. 

If $G={\rm S}_4$ then the involution $y$ must be a transposition: the double transpositions lie in the normal subgroup $K\cong {\rm V}_4$ 
with $G/K\cong{\rm S}_3$ non-cyclic, so they cannot be members of generating pairs. One can take $x$ to be any $3$-cycle or $4$-cycle not 
inverted by $y$, giving the cube $\{4,3\}$ or its dual, the octahedron $\{3,4\}$. The operations $H_j$ act trivially on $\O(G)$.

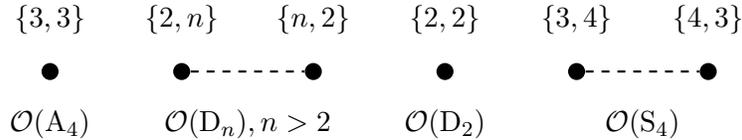
\begin{figure}[h!]
\begin{center}
\begin{tikzpicture}[scale=0.35, inner sep=0.8mm]

\node (a) at (0,0) [shape=circle, fill=black] {};
\node (b) at (5,0) [shape=circle, fill=black] {};
\node (c) at (10,0) [shape=circle, fill=black] {};
\node (d) at (15,0) [shape=circle, fill=black] {};
\node (e) at (20,0) [shape=circle, fill=black] {};
\node (f) at (25,0) [shape=circle, fill=black] {};

\draw [thick, dashed] (b) to (c);
\draw [thick, dashed] (e) to (f);

\node at (0,-2) {$\O({\rm A}_4)$};
\node at (7.5,-2) {$\O({\rm D}_n), n>2$};
\node at (15,-2) {$\O({\rm D}_2)$};
\node at (22.5,-2) {$\O({\rm S}_4)$};

\node at (0,2) {$\{3,3\}$};
\node at (5,2) {$\{2,n\}$};
\node at (10,2) {$\{n,2\}$};
\node at (15,2) {$\{2,2\}$};
\node at (20,2) {$\{3,4\}$};
\node at (25,2) {$\{4,3\}$};

\end{tikzpicture}

\end{center}
\caption{The graphs $\O({\rm A}_4)$,  $\O({\rm D}_n)$ and $\O({\rm S}_4)$}
\label{fig:A4}
\end{figure}

The graphs $\O(G)$ corresponding to these groups $G$ are shown in Figure~\ref{fig:A4}.

\begin{figure}[h!]
\begin{center}
\begin{tikzpicture}[scale=0.3, inner sep=0.8mm]

\node (a) at (5,0) [shape=circle, fill=black] {};
\node (b) at (0,0) [shape=circle, fill=black] {};
\node (c) at (-5,0) [shape=circle, fill=black] {};

\draw [thick] (a) to (b);
\draw [thick, dashed] (b) to (c);

\node at (5,2) {$\{5,5\}_6$};
\node at (0,2) {$\{3,5\}$};
\node at (-5,2) {$\{5,3\}$};

\end{tikzpicture}

\end{center}
\caption{The graph $\O({\rm A}_5)$} 
\label{fig:A5}
\end{figure}
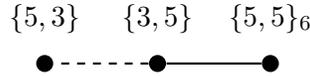

If $G$ is the alternating group ${\rm A}_5$ (also isomorphic to $\PSL_2(4)$ and $\PSL_2(5)$) then $\O(G)$ consists of three maps, 
namely the dodecahedron $\{5,3\}$, the icosahedron $\{3,5\}$ and the great dodecahedron $\{5,5\}_6$ of genus~$4$ (the self-dual 
map R4.6 in~\cite{Con09}); duality $D$ transposes the first two, while $H_2$ transposes the last two, so $\O(G)$ is connected (see 
Figure~\ref{fig:A5}).

\section{Isotactic polygons}\label{sec:isotactic}

Let $\M$ be an orientably regular map of type $\{p,q\}$. An {\em isotactic polygon\/} $P$ of {\em type\/} $(d_1,d_2, \ldots, d_m)^n$ in $\M$ is 
a closed path in the underlying graph of $\M$ formed by successively taking the $d_i$-th edge at the $i$-th vertex $v_i$ visited, following the local 
orientation around $v_i$, using the sequence $d_1, d_2, \ldots, d_m\in{\mathbb Z}_q$ repeated $n$ times. We call $d_i$ the {\em right degree\/} 
of $P$ at $v_i$, and $l=mn$ the {\em length\/} of $P$. By orientable regularity, if such a polygon $P$ exists, then it does so starting at any directed 
edge in $\M$, since it is equivalent to a relation
\[(x^{d_1}yx^{d_2}y\ldots x^{d_m}y)^n=1\]
in the monodromy group $G=\langle x, y\rangle$ of $\M$.

This concept is a common generalisation of the classical notions of {\em face\/}, {\em Petrie polygon\/} and {\em hole\/}. In fact, 
\begin{itemize}
\item an isotactic polygon of type $(1)^p$ or $(-1)^p$ is the boundary of a $p$-valent face;
\item more generally, an isotactic polygon of type $(j)^l$ with $j$ coprime to $q$ is a $j$-hole of length $l$;
\item an isotactic polygon of type $(1,-1)^{r/2}$ is a Petrie polygon of length $r$ (provided this length is even);
\item more generally, an isotactic polygon of type $(j,-j)^{r/2}$ is a Petrie polygon of order $j$, where $r$ is its (even) length.
\end{itemize}

\begin{figure}[!h]
\begin{center}
\includegraphics[width=0.65\textwidth] {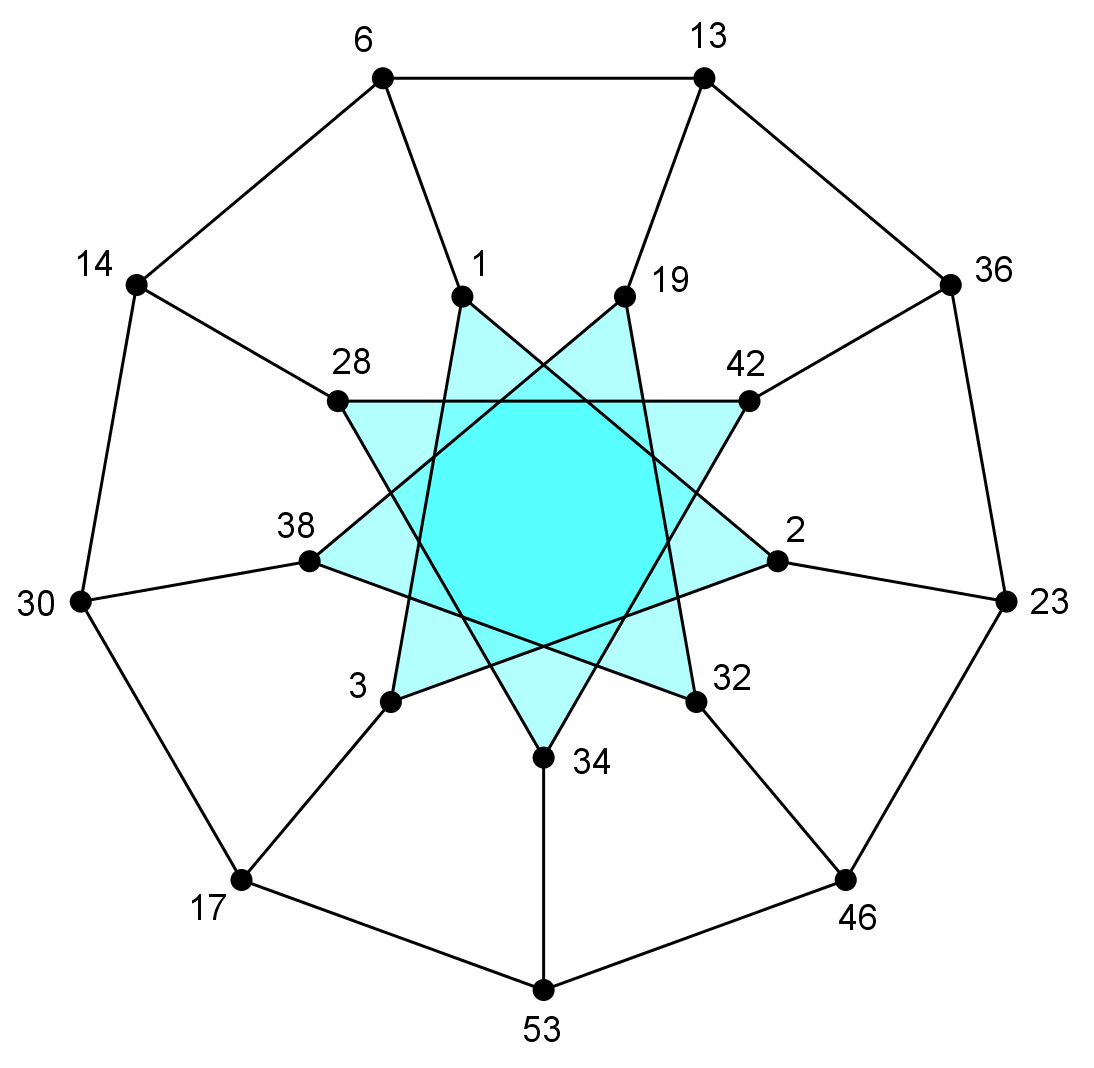}
\caption{A generalized Petersen graph of type GP(9,3) in the Fricke--Macbeath map $\Fr$ of type $\{3,7\}$ and genus 7 (vertex labels correspond to those given in Figure~\ref{fig:H7} below).}
\label{fig:GP}
\end{center}
\end{figure}

These well-known examples of isotactic polygons are distinguished by their roles in the Petrie and hole operations, but other types also occur in various roles. 
For example, hexagons of type $(2,4)^3$ occur in the polyhedral realisation of Klein's map $\K$ of genus $3$ (see Section~\ref{sec:PSL_2(7)}) given by 
Schulte and Wills in~\cite{SW85}. On the other hand, the heuristic role of hexagons of type $(2,4,4)^2$ is emphasised in the investigation of $7$-fold 
rotational symmetry of the Fricke--Macbeath map $\Fr$ of genus $7$ (see Section~\ref{sec:PSL_2(8)}) in~\cite{BC17, BG}. Moreover, in this map the 
$3$-holes, together with three suitable triangular faces, form generalised Petersen graphs of type ${\rm GP}(9,3)$ (see Figure~\ref{fig:GP}); the presence 
of such subgraphs in the underlying graph of $\Fr$ confirms that this map has no polyhedral embedding in ${\mathbb E}^3$ with $9$-fold symmetry~\cite{BG}.

\section{Regularity}\label{sec:regular}

An orientably regular map $\M$, corresponding to a generating pair $x, y$ of $G$, is regular if and only if there is an automorphism $\alpha$ of $G$ inverting $x$ and $y$, 
or equivalently inverting $x$ and centralising $y$. In this case the full automorphism group ${\rm Aut}\,\M$ is a semidirect product $G\rtimes\langle\alpha\rangle$ of $G$ 
and ${\rm C}_2$, with $\alpha$ acting a reflection of $\M$. For example, the maps in $\O(G)$ for $G={\rm A}_4$, ${\rm D}_n$ and ${\rm S}_4$ are all regular, with 
${\rm Aut}\,\M\cong{\rm S}_4$, ${\rm D}_n\times{\rm C}_2$ and ${\rm S}_4\times{\rm C}_2$ respectively.

We will say that a regular map $\M$ is {\em inner\/} or {\em outer regular\/} as $\alpha$ is or is not an inner automorphism, that is, conjugation by some element 
$c\in G$. If $\alpha$ is inner then $c^2$ is in the centre $Z(G)$ of $G$, so if $Z(G)=1$ (as is the case with all the groups considered here, apart from ${\rm D}_n$ 
for $n\le 2$), then $c^2=1$ and ${\rm Aut}\,\M\cong G\times{\rm C}_2$ with $\alpha c=\alpha c^{-1}$ generating the second direct factor. In this case $\M$ has 
a non-orientable regular quotient map $\overline\M=\M/{\rm C}_2$, of genus $g+1$ where $\M$ has genus $g$, with automorphism group $G$; then $\M$ is the 
canonical orientable double cover of $\overline\M$.

It is easy to check that regularity, inner regularity and ${\rm Aut}\,\M$ are preserved by the operations $D$ and $H_j$, so they are constant throughout each 
connected component of $\O(G)$. In particular, maps in $\O(G)$ with different regularity properties must lie in different connected components of the graph.


\subsection{${\mathcal O}(G)$ for $G={\rm S}_5$}\label{subsec:S5}

The earlier examples were very straightforward, involving familiar maps. This case needs rather more thought.

If $G={\rm S}_5$ then the involution $y$ can be a transposition or a double transposition. In the first case, there are possible generators $x$ of orders $4$, $5$ and $6$, 
giving rise to two dual pairs of regular maps, the dual pair R4.2 of types $\{4,5\}_6$ and $\{5,4\}_6$ on Bring's curve of genus~$4$ (see~\cite{Web05}, for example), 
and the dual pair R9.16 of types $\{5,6\}_4$ and $\{6,5\}_4$, with $H_2$ transposing R4.2a and R9.16b; the component of $\O(G)$ containing these four maps is a 
path graph, with three edges labelled $D, H_2, D$.

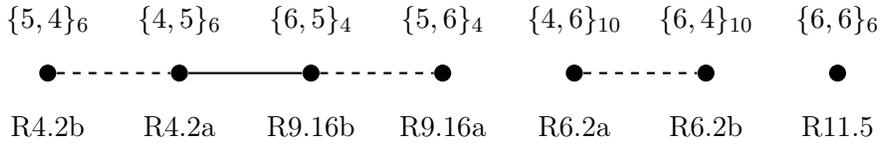
\begin{figure}[h!]
\begin{center}
\begin{tikzpicture}[scale=0.35, inner sep=0.8mm]

\node (a) at (0,0) [shape=circle, fill=black] {};
\node (b) at (5,0) [shape=circle, fill=black] {};
\node (c) at (10,0) [shape=circle, fill=black] {};
\node (d) at (15,0) [shape=circle, fill=black] {};
\node (e) at (20,0) [shape=circle, fill=black] {};
\node (f) at (25,0) [shape=circle, fill=black] {};
\node (g) at (30,0) [shape=circle, fill=black] {};

\draw [thick, dashed] (a) to (b);
\draw [thick] (b) to (c);
\draw [thick, dashed] (c) to (d);
\draw [thick, dashed] (e) to (f);

\node at (0,-2) {R4.2b};
\node at (5,-2) {R4.2a};
\node at (10,-2) {R9.16b};
\node at (15,-2) {R9.16a};
\node at (20,-2) {R6.2a};
\node at (25,-2) {R6.2b};
\node at (30,-2) {R11.5};

\node at (0,2) {$\{5,4\}_6$};
\node at (5,2) {$\{4,5\}_6$};
\node at (10,2) {$\{6,5\}_4$};
\node at (15,2) {$\{5,6\}_4$};
\node at (20,2) {$\{4,6\}_{10}$};
\node at (25,2) {$\{6,4\}_{10}$};
\node at (30,2) {$\{6,6\}_6$};

\end{tikzpicture}
\end{center}
\caption{The graph $\O({\rm S}_5)$} 
\label{fig:S5}
\end{figure}

In the second case $x$, which must be odd, can have order $4$ or $6$, giving the dual pair R6.2 of types $\{4,6\}_{10}$ and $\{6,4\}_{10}$,  together with the 
self-dual map R11.5 of type $\{6,6\}_6$ (for this last map one can take $x=(1,5,3)(2,4), y=(1,2)(3,4)$); all operations $H_j$ act trivially on these maps, so the 
pair R6.2 and the map R11.5 form two more connected components of $\O(G)$. The graph is shown in Figure~\ref{fig:S5}: the vertices are labelled with the extended 
types and our extension of the notation in~\cite{Con09} for the corresponding maps. This example shows that the involution $y$ (or more precisely its orbit under 
${\rm Aut}\,G$) is insufficient to characterise a component of $\O(G)$.

The maps $\M\in\O(G)$ are all regular, and since ${\rm Out}\,{\rm S}_5$ and $Z({\rm S}_5)$ are both trivial they are inner regular, with automorphism group ${\rm S}_5\times{\rm C}_2$;
their non-orientable quotients $\overline\M$, reading Figure~\ref{fig:S5} from left to right, are the dual pairs N5.1, N10.4, N7.1 and the self-dual map N12.3, all with automorphism group 
${\rm S}_5$.

\medskip

\begin{figure}[!h]
\begin{center}
\includegraphics[width=0.85\textwidth] {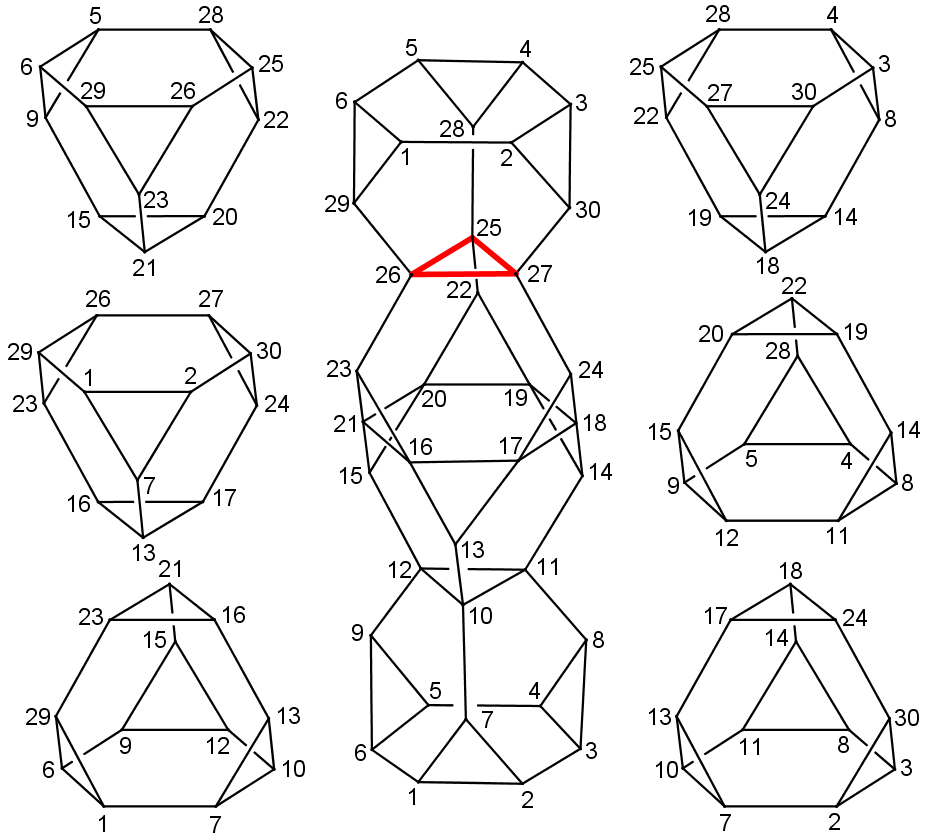}
\caption{The combinatorial structure of R6.2b~\cite{GG}, with one hole of length 3 highlighted.}
\label{fig:tencell}
\end{center}
\end{figure}
We note that R6.2b occurs first as one of Coxeter's regular skew polyhedra~\cite{Cox}. It is realised as a subcomplex of the boundary complex of the dual 
of the 4-polytope $t_{1,2}\{3,3,3\}$. Its faces form the faces of 10 equal Archimedean truncated tetrahedra (the facets of the polytope). The 
non-identity element of ${\rm C}_2$ is realised in this case as a central inversion (i.e.\ reflection in a point) in $\mathbb E^4$ interchanging 
two disjoint 10-tuples of hexagons. Two triples of hexagons sharing a hole of length 3 belong to the same such 10-tuples; in Figure~\ref{fig:tencell} 
both the hexagonal faces and the triangular holes can clearly be seen.

R4.2a goes back even earlier, in fact to Gordan's work~\cite{Gor}. Both the dual pairs R4.2a, R4.2b and R6.2a, R6.2b have polyhedral realisations 
in Euclidean 3-space~\cite{SW86, SW87}.


\section{Hurwitz surfaces, groups and maps}\label{sec:Hurwitz}

One aim of this paper is to study $\O(G)$ where $G$ is a Hurwitz group, so here we summarise some of the important properties of this class of groups.

Hurwitz~\cite{Hur93} showed that the automorphism group $G={\rm Aut}\,\S$ of a compact Riemann surface $\S$ of genus $g\ge 2$ has order at most $84(g-1)$, attained 
if and only if $\S\cong \H/K$, where $\H$ is the hyperbolic plane and $K$ is a normal subgroup of finite index in the triangle group
\[\Delta=\Delta(7,2,3)=\langle X, Y, Z\mid X^7=Y^2=Z^3=XYZ=1\rangle\]
with $G\cong \Delta/K$. Such surfaces and groups attaining this upper bound are called {\sl Hurwitz surfaces\/} and {\sl Hurwitz groups}. In each case, $\S$ carries an orientably 
regular map $\M$ of type $\{3,7\}$, called a {\sl Hurwitz map}, with orientation-preserving automorphism group $G$. This map is regular if and only if $K$ is normal in the extended 
triangle group $\Delta[7,2,3]$ which contains $\Delta$ with index $2$, in which case the full automorphism group of $\M$ is isomorphic to $\Delta[7,2,3]/K$; this is equivalent to $G$ 
having an automorphism inverting two of its standard generators $x$, $y$ and $z$. Conder has written some very useful surveys on Hurwitz groups in~\cite{Con90, Con10}.

Every Hurwitz group is perfect, since $\Delta$ is, so it is a covering of a non-abelian finite simple group, which is itself a Hurwitz group. Two important infinite classes of simple 
Hurwitz groups are given by the following theorems~\cite{Macb69, Con80}:

\begin{theorem}[Macbeath, 1969]\label{th:Macb}
The group $\PSL_2(q)$ is a Hurwitz group if and only if
\begin{enumerate}
\item $q=7$, with a unique Hurwitz surface and map, or
\item $q=p$ for some prime $p\equiv\pm 1$ {\rm mod}~$(7)$, with three Hurwitz surfaces and maps for each $q$, or
\item $q=p^3$ for some prime $p\equiv\pm 2$ or $\pm 3$ {\rm mod}~$(7)$, with one Hurwitz surface and map for each $q$.
\end{enumerate}
\end{theorem}

\begin{theorem}[Conder, 1980]\label{th:Con}
The alternating group ${\rm A}_n$ is a Hurwitz group for each $n\ge 168$.
\end{theorem}

In~\cite{Con80}, Conder also determined which alternating groups of degree $n<168$ are Hurwitz groups; the smallest is ${\rm A}_{15}$.


\section{Maps with $G=\PSL_2(q)$}\label{sec:PSL_2(q)}

\subsection{Properties of ${\rm PSL}_2(q)$}

Later in this paper we will construct the graphs ${\mathcal O}(G)$ for some specific groups $G={\rm PSL}_2(q)$. Here we briefly summarise 
a few of the properties of these groups which we will need; see~\cite[Ch.~XII]{Dic} or~\cite[II.8]{Hup} for full details. 

The group $G=\PSL_2(q)=\SL_2(q)/\{\pm I\}$, $q=p^e$, has order $q(q^2-1)$ or $q(q^2-1)/2$ as $p=2$ or $p>2$. It is simple for all $q\ge 4$.

One useful property of $G=\PSL_2(q)$ is that conjugacy classes of non-identity elements $g\in G$ are determined uniquely by their traces 
(more precisely, trace-pairs) $\pm{\rm tr}(g)\in \F_q$. For example, the elements of orders $2$, $3$ and $p$ are the non-identity elements 
with traces $0$, $\pm 1$ and $\pm 2$ respectively.


\subsection{Consistent choice of $y$}

For any group $G$, since the generator $y$ is invariant (up to automorphisms) under the operations $D$ and $H_j$, it makes sense to use the same 
element for each map in a given component of $\O(G)$, or even for each map in $\O(G)$ when this is possible.

For example, if $G={\rm PSL}_2(q)$ then all involutions in $G$ are conjugate to
\[y=\begin{pmatrix}0 & 1 \\ -1 & 0\end{pmatrix}.\]
(Here we identify each matrix $A\in{\rm SL}_2(q)$ with $-A$.)
With this as $y$, its zero entries make it easier to see the traces of words such as $x^jy$ and $[x,y]$, and hence to  find their orders, giving the extended types of various maps.

For instance, if $y$ is as above and
\[x=\begin{pmatrix}a & b \\ c & d\end{pmatrix}\quad\hbox{then}\quad xy=\begin{pmatrix}-b & a \\ -d & c\end{pmatrix},
\quad\hbox{so}\quad z=(xy)^{-1}=\begin{pmatrix}c & -a \\ d & -b\end{pmatrix},\]
with trace $\pm(c-b)$. We will call $\tau=\pm(a+d)$ and $\tau'=\pm(c-b)$ the {\sl trace\/} and {\sl cotrace\/} of the element $x$ and of the triple $(x,y,z)$. 
If this triple generates $G$ then it corresponds to a map $\M\in\O(G)$ of type $\{p,q\}$, where the pair consisting of the trace $\pm(a+d)$ and cotrace 
$\pm(c-b)$ (unique up to the action of ${\rm Gal}\,\F_q$) determine the orders $q$ and $p$ of $x$ and $z$ and thus the type $\{p,q\}$ of $\M$, together 
with its genus
\[g=1+\frac{|G|}{2(\frac{1}{2}-\frac{1}{p}-\frac{1}{q})}\]
by the Riemann--Hurwitz formula. We also have
\[[x,y]=x^{-1}y^{-1}xy=x^{-1}yxy=\begin{pmatrix}-b^2-d^2 & ab+cd \\ ab+cd & -a^2-c^2\end{pmatrix}\]
with trace $\sigma=\pm(a^2+b^2+c^2+d^2)$, giving the Petrie length $r$ (twice the order of $[x,y]$) and hence the extended type $\{p,q\}_r$ of $\M$.

Once $y$ is chosen, for example as above, all of this data for a map $\M$ is uniquely determined by the generating element $x$, so instead of working with 
triples $(x,y,z)$ it is more efficient to work with the single elements $x$ when applying operations such as $D$ and $H_j$.

For example, if $\tau=\pm(a+d)$ is the trace of $x$ then
\[x^2=\begin{pmatrix}a^2+bc & b(a+d) \\ c(a+d) & d^2+bc\end{pmatrix}\]
has trace $\pm(a^2+2bc+d^2)=\pm(a^2+2(ad-1)+d^2)=\pm(\tau^2-2)$, and
\[(x^2y)^{-1}=\begin{pmatrix}c(a+d) & -a^2-bc \\ d^2+bc & -b(a+d)\end{pmatrix},\]
has trace $\pm(a+d)(c-b)=\pm\tau\tau'$, giving the cotrace and hence the type and genus of $H_2(\M)$.
Similarly, its Petrie length is determined by the trace
\[\pm((a^2+bc)^2+(b(a+d))^2+(c(a+d))^2+(d^2+bc)^2)\]
of $[x^2,y]$. Applying $H_{-1}$ or $D$ is achieved by replacing $x$ with
\[x^{-1}=\begin{pmatrix}d & -b \\ -c & a\end{pmatrix} \quad\hbox{or}\quad
z=\begin{pmatrix}c & -a \\ d & -b\end{pmatrix}.\]
Of course, these leave $r$ and $g$ invariant, while $D$ transposes $p$ and $q$.

By an observation of Singerman~\cite{Sin74}, if $G={\rm PSL}_2(q)$ then all maps $\M\in\O(G)$ are regular, so $H_{-1}(\M)\cong \M$ and there is 
no need to apply $H_{-1}$. However, this allows us, if we wish, to replace $x$ with
\[z^{-1}=xy=\begin{pmatrix}-b & a \\ -d & c\end{pmatrix}\]
instead of $z$ when applying $D$.

Of course, many of the above expressions for matrices and traces are a little simpler if $q$ is a power of $2$.


\section{$G={\rm PSL}_2(7)$}\label{sec:PSL_2(7)}

There is a unique Hurwitz surface of genus~$3$, Klein's quartic curve (see~\cite{Kle78, Lev}) given in projective coordinates by $x_0^3x_1+x_1^3x_2+x_2^3x_0=0$ . 
The corresponding Hurwitz map $\K$ of type $\{3,7\}$, {\sl Klein's map} (see Figure~\ref{fig:Klein}), has orientation-preserving automorphism group 
$G\cong {\rm PSL}_2(7)\cong {\rm GL}_3(2)$, a simple group of order $168$, and has full automorphism group ${\rm PGL}_2(7)$. This is the map R7.1 
in~\cite{Con09}, discussed as a polyhedron by Schulte and Wills in~\cite{SW85}. 

The non-identity elements of ${\rm PSL}_2(7)$ have trace $0$, $\pm 1$, $\pm 2$ or $\pm3$ as they have order $2$, $3$, $7$ or $4$ respectively. 
For a Hurwitz triple we require $x$ and $z$ to have orders $7$ and $3$, that is, $\tau=\pm(a+d)=\pm 2$, $\tau'=\pm(c-b)=\pm 1$ and of course 
$ad-bc=1$, so we may take $a=1$, $b=1$, $c=0$, $d=1$; the resulting generating triple for $\K$ is
\[x=\begin{pmatrix}1 & 1 \\ 0 & 1\end{pmatrix},\quad y=\begin{pmatrix}0 & 1 \\ -1 & 0\end{pmatrix},
\quad z=\begin{pmatrix}0 & -1 \\ 1 & -1\end{pmatrix}.\]
Since $\sigma=a^2+b^2+c^2+d^2=3$, $[x,y]$ has order $4$, so $\K$ has Petrie length $8$ and extended type $\{3,7\}_8$.

Squaring $x$ to apply $H_2$, we obtain
\[x^2=\begin{pmatrix}1 & 2 \\ 0 & 1\end{pmatrix},\]
corresponding to the map $H_2(\K)$, which has extended type $\{7,7\}_6$ and genus~$19$. Squaring again, we obtain
\[x^4=\begin{pmatrix}1 & 4 \\ 0 & 1\end{pmatrix},\]
corresponding to the map $H_4(\K)=H_{-3}(\K)=H_3(\K)$ of type $\{4,7\}_8$ and genus $10$. Since $2^3\equiv 1$ mod~$(7)$, applying $H_2$ again simply returns us to $\K$.

For $D(\K)$, of type $\{7,3\}_8$ and genus $3$, we can replace $x$ with
\[z=\begin{pmatrix}0 & -1 \\ 1 & -1\end{pmatrix},\]
and for $D(H_4(\K))$, of type $\{7,4\}_8$ and genus $10$, we can replace $x^4$ with
\[\begin{pmatrix}0 & -1 \\ 1 & 3\end{pmatrix}.\]
However, $H_2(\K)$ is self-dual, that is, $DH_2(\K)\cong H_2(\K)$; this follows immediately from the following result:
 
\begin{proposition}\label{prop:PSL_2(7)}
There is a unique orientably regular map of type $\{7,7\}$ with orientation-preserving automorphism group $G\cong{\rm PSL}_2(7)$.
\end{proposition}

\noindent{\sl Proof.}  Any such map $\M$ corresponds to a generating triple $(x,y,z)$ for $G$ of type $(7,2,7)$. We will use the Frobenius triple-counting 
formula~\cite{Fro96}, which states that if $\mathcal X$, $\mathcal Y$ and $\mathcal Z$ are conjugacy classes in any finite group $G$, then the number 
of triples $x\in{\mathcal X}$, $y\in{\mathcal Y}$ and $z\in{\mathcal Z}$ with $xyz=1$ is equal to
 \begin{equation}\label{eq:Frobenius}
 \frac{|{\mathcal X}|\cdot |{\mathcal Y}| \cdot |{\mathcal Z}|}{|G|}\sum_{\chi}\frac{\chi(x)\chi(y)\chi(z)}{\chi(1)},
 \end{equation}
 where the sum is over all irreducible complex characters $\chi$ of $G$.
 
 In ATLAS notation~\cite{ATLAS}, the conjugacy classes $\mathcal X$, $\mathcal Y$ and $\mathcal Z$ of $G={\rm PSL}_2(7)$ containing $x$, $y$ and $z$ 
are respectively 7A or 7B, 2A, and 7A or 7B. In all cases
 \[\frac{|{\mathcal X}| \cdot |{\mathcal Y}| \cdot |{\mathcal Z}|}{|G|}
 =\frac{(2^3\cdot 3)\cdot(3\cdot 7)\cdot(2^3\cdot 3)}{2^3\cdot 3\cdot 7}
 =2^3\cdot 3^2 = 72.\]
 If ${\mathcal X}\ne{\mathcal Z}$ then by using the character values in~\cite{ATLAS} we find that the character sum in (\ref{eq:Frobenius}) is
 \[1+2\cdot\frac{(-1)\cdot (-1+i\sqrt{7})/2 \cdot (-1-i\sqrt{7})/2}{3}+\frac{(-1)\cdot 2\cdot (-1)}{6}=0.\]
Thus there are no such triples in $G$ and hence there are no corresponding maps. If ${\mathcal X}={\mathcal Z}=$ 7A or 7B then the character sum is
 \[1+\frac{(-1)\cdot ((-1+i\sqrt{7})/2)^2}{3}+\frac{(-1)\cdot ((-1-i\sqrt{7})/2)^2}{3}+\frac{(-1)\cdot 2\cdot (-1)}{6}=\frac{7}{3},\]
 so, adding these contributions, the total number of such triples is
 \[2\cdot 72\cdot \frac{7}{3}=336.\]
 
These triples all generate $G$, since no proper subgroup of $G$ has order divisible by $14$. The triples are permuted by ${\rm Aut}\,G={\rm PGL}_2(7)$, 
semi-regularly since only the identity automorphism can fix a generating set. Since $|{\rm PGL}_2(7)|=336$, ${\rm Aut}\,G$ is transitive on the triples, 
so $\Delta(7,2,7)$ has one normal subgroup with quotient $G$, and hence there is one orientably regular map of type $\{7,7\}$ with orientation-preserving 
automorphism group $G$. \hfill$\square$

\medskip

\begin{figure}[!h]
\begin{center}
\includegraphics[width=0.85\textwidth] {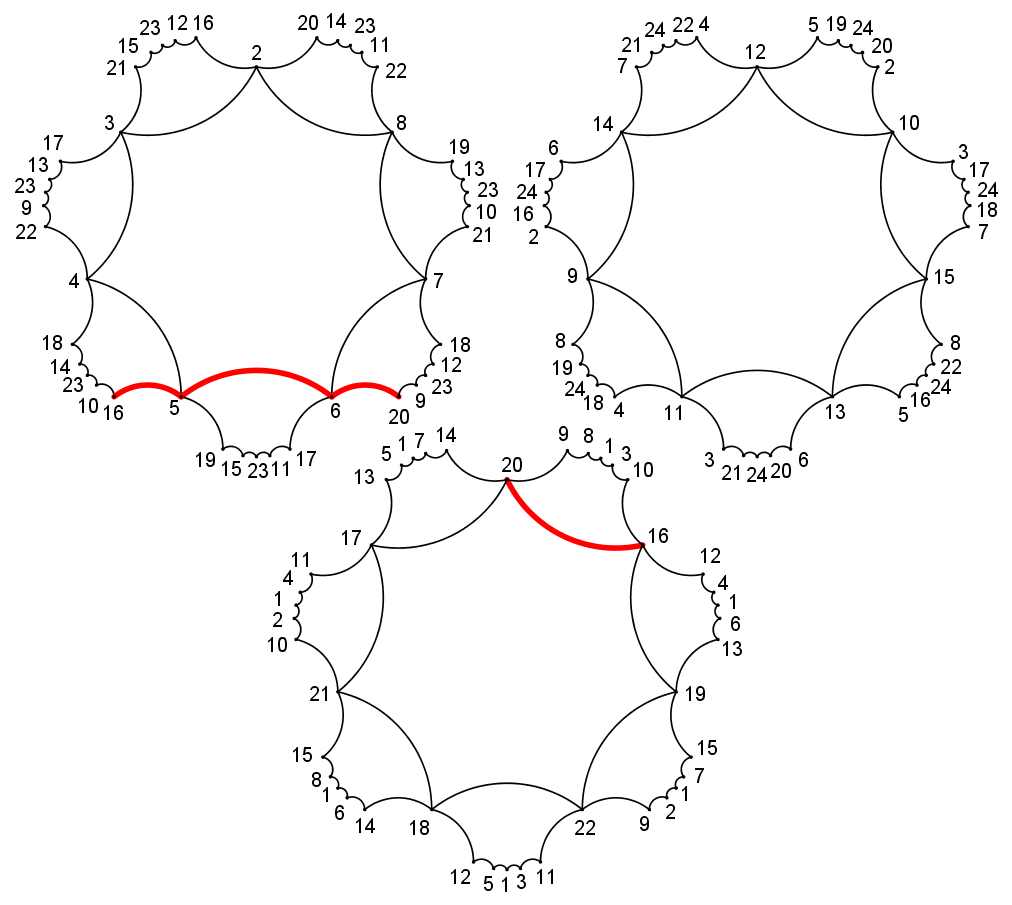}
\caption{The map $H_2(\K)$ of type $\{7,7\}_6$ and genus $19$. It is trisected, and
each part is represented in the disc model of the hyperbolic plane. A
2-hole is highlighted; this is the same as the 3-hole in $\K$ highlighted in Figure~\ref{fig:Klein}.}
\label{fig:77}
\end{center}
\end{figure}

By the Riemann--Hurwitz formula, the above map $\M$ has genus $19$. By its uniqueness, it is regular and self-dual, and is isomorphic to $H_2(\K)$ 
and to the map R19.23 of extended type $\{7,7\}_6$, the only regular orientable map of this genus and type in~\cite{Con09}. This map is shown in 
Figure~\ref{fig:77}.

\medskip

\medskip
 
Our other claims for uniqueness of maps, or enumerations of maps of a given genus and type, can be proved in a similar way. However, in future we will not give 
such full details unless there are special aspects of the calculation which need to be mentioned.

In fact the maps $\K$, $D(\K)$, $H_2(\K)$, $H_4(K)$ and $DH_4(\K)$ are the only maps in $\O(G)$. To see this, note that the face- and vertex-valencies $p$ 
and $q$ for any such map must be orders of non-identity elements of $G$, so they must take the values $2, 3, 4$ or $7$. Since $p^{-1}+q^{-1}<1/2$ the only 
possible types (up to duality) are $\{3,7\}$, $\{7,3\}$, $\{4,7\}$, $\{7,4\}$ and $\{7,7\}$; as in the case of $\{7,7\}$ it can be verified by using the Frobenius 
triple-counting formula (or by inspection of~\cite{Con09}) that there is only one map of each type in $\O(G)$.

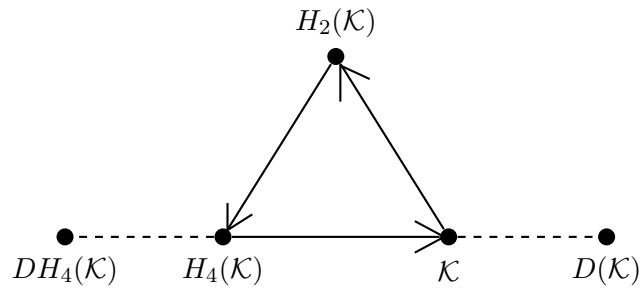
\begin{figure}[h!]
\begin{center}
\begin{tikzpicture}[scale=0.3, inner sep=0.8mm]

\node (a) at (5,0) [shape=circle, fill=black] {};
\node (b) at (0,8) [shape=circle, fill=black] {};
\node (c) at (-5,0) [shape=circle, fill=black] {};
\node (d) at (12,0) [shape=circle, fill=black] {};
\node (e) at (-12,0) [shape=circle, fill=black] {};
\draw [thick] (a) to (b) to (c) to (a);
\draw [thick, dashed] (a) to (d);
\draw [thick, dashed] (c) to (e);

\draw [thick] (0.2,6) to (0.2,7.6) to (1.5,7);
\draw [thick] (-4.8,1.5) to (-4.8,0.3) to (-3.7,1);
\draw [thick] (3.5,0.7) to (4.8,0) to (3.5,-0.7);

\node at (5,-1.5) {$\K$};
\node at (12,-1.5) {$D(\K)$};
\node at (0,9.5) {$H_2(\K)$};
\node at (-5,-1.5) {$H_4(\K)$};
\node at (-12,-1.5) {$DH_4(\K)$};

\end{tikzpicture}

\end{center}
\caption{The graph $\O({\rm PSL}_2(7))$} 
\label{fig:PSL_2(7)}
\end{figure}

The graph $\O(G)$ is therefore as shown in Figure~\ref{fig:PSL_2(7)}, where the actions of $D$ and $H_2$ are represented by undirected broken and directed 
unbroken edges; that of $H_3$ is given by reversing the directed edges for $H_2$.

In Table~\ref{OPSL27} we list the maps $\M\in\O({\rm PSL}_2(7))$. The first column shows the corresponding entry in~\cite{Con09}, with letters `a' and `b' assigned 
as explained earlier. The second column gives the extended type $\{p,q\}_r$ of $\M$, with $p$, $q$ and $r$ denoting the face- and vertex-valencies and the Petrie length. 
The third column shows how $\M$ may be obtained from Klein's map $\K$ by applying duality and hole operations. The final column gives the effect of applying the hole 
operation $H_2$ to $\M$; this is left blank if $\M$ is unchanged (as when $q=3$, for example) or $q$ is even. Of course, duality $D$ transposes pairs R$g.n$a and R$g.n$b, 
while leaving maps R$g.n$ invariant. Since the graph $\O(G)$ is connected, and $\K$ is regular with full automorphism group ${\rm PGL}_2(7)$, the same applies to every 
map $\M$ in $\O(G)$.
 
 \begin{table}[ht]
\centering
\begin{tabular}{| p{2cm} | p{1.3cm} | p{3cm} | p{1.2cm} |}
\hline
Entry in~\cite{Con09} & Type & Relationship to $\K$ & $H_2(\M)$ \\
\hline\hline
R3.1a & $\{3,7\}_8$ & $\K$ & R19.23\\
\hline
R3.1b & $\{7,3\}_8$ & $D(\K)$ &  \\
\hline
R10.9a & $\{4,7\}_8$ & $H_4(\K)$ & R3.1a \\
\hline
R10.9b & $\{7,4\}_8$ & $DH_4(\K)$ & \\
\hline
R19.23 & $\{7,7\}_6$ & $H_2(\K)$ & R10.9a \\
\hline

\end{tabular}
\caption{The maps $\M$ in $\O({\rm PSL}_2(7))$.}
\label{OPSL27}
\end{table}

\medskip

\noindent{\bf Example} Taking subscripts mod~$(7)$ and using regularity, we see that
\[H_3(\K)\cong H_{-3}(\K)\cong H_4(\K)\cong H_2^2(\K)\cong{\rm R10.9a},\]
of type $\{4,7\}$; Figure~\ref{fig:Klein} confirms this, showing that the 3-holes of $\K$, giving the faces of $H_3(\K)$, have length $4$.


\section{$G={\rm PSL}_2(8)={\rm SL}_2(8)$}\label{sec:PSL_2(8)}

\begin{figure}[!h]
\begin{center}
\includegraphics[width=0.975\textwidth] {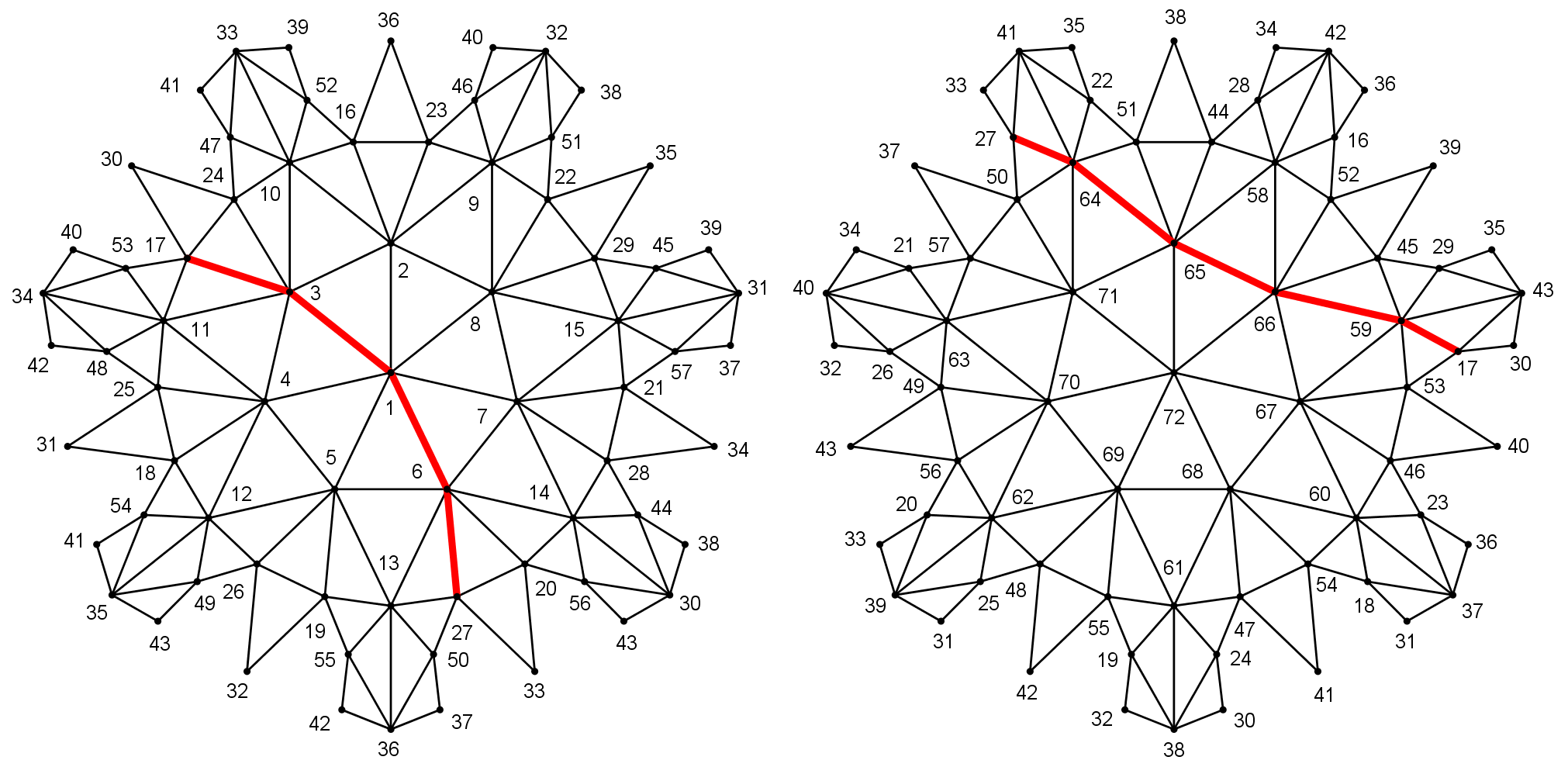}
\caption{Combinatorial scheme of the Hurwitz map $\Fr$ of genus $7$ (with vertex labels taken from \cite{BG}).
A 3-hole is highlighted.}
\label{fig:H7}
\end{center}
\end{figure}

There is a unique Hurwitz surface $\S$ of genus~$7$, described by Fricke~\cite{Fri99} in 1899 and rediscovered by Macbeath~\cite{Macb65} 
in 1965; its automorphism group $G$, isomorphic to the simple group ${\rm PSL}_2(8)={\rm SL}_2(8)$ of order $504$, is the orientation-preserving 
automorphism group of a regular map $\Fr$ of type $\{3,7\}$ on $\S$, the {\sl Hurwitz map\/} of genus $7$, with full automorphism group 
$G\times{\rm C}_2$. The combinatorial structure of this map is shown in Figure~\ref{fig:H7}. For its polyhedral realisations, 
see~\cite{BC17, BG, Bok21}.

\subsection{The group $G$}

In order to work with this group $G$, let us represent the field $\F_8$ of order $8$ as $\F_2[t]/(t^3+t+1)$, with its elements represented as polynomials 
in $\F_2[t]$ of degree at most $2$. Then the multiplicative group of $\F_8$ consists of
\[t,\; t^2,\; t^3=t+1,\; t^4=t^2+t,\; t^5=t^2+t+1,\; t^6=t^2+1,\; t^7=1.\]
The non-identity elements of $G$ form the following conjugacy classes:
\begin{itemize}
\item one class of elements of order $2$, with trace $0$;
\item one class of elements of order $3$, with trace $1$;
\item three classes of elements of order $7$, with traces $t+1, t^2+1, t^2+t+1$;
\item three classes of elements of order $9$, with traces $t, t^2, t^2+t$.
\end{itemize}
Each class is inverse-closed. The outer automorphism group of $G$, isomorphic to ${\rm C}_3$, is induced by the Galois group of $\F_8$ which is generated by 
the Frobenius automorphism $t\mapsto t^2$. This group permutes the three classes of elements of order $7$ in a single cycle, and likewise  for those of order $9$.

Since ${\rm Out}\,G$ has odd order, and $Z(G)$ is trivial, all maps in $\O(G)$ are inner regular, with full automorphism group $G\times{\rm C}_2$.


\subsection{Construction of $\O(G)$}

Since all involutions in $G$ are conjugate, in forming $\O(G)$ we can take
\[y=\begin{pmatrix}0 & 1 \\ -1 & 0\end{pmatrix}=\begin{pmatrix}0 & 1 \\ 1 & 0\end{pmatrix},\]
using the fact that $\F_8$ has characteristic $2$ to eliminate minus signs. If
\[x=\begin{pmatrix}a & b \\ c & d\end{pmatrix}\quad\hbox{then}\quad z=\begin{pmatrix}c & a \\ d & b\end{pmatrix},\]
with the trace $\tau=a+d$ and the cotrace $\tau'=b+c$ giving the type and hence the genus of $\M$. Also the trace $\sigma=a^2+b^2+c^2+d^2$ of $[x,y]$ 
gives its Petrie length and hence its extended type. Then
\[x^2=\begin{pmatrix}a^2+bc & b(a+d) \\ c(a+d) & d^2+bc\end{pmatrix}\]
has trace $a^2+d^2=\tau^2$, and the trace of 
\[(x^2y)^{-1}=\begin{pmatrix}c(a+d) & a^2+bc \\ d^2+bc & b(a+d)\end{pmatrix}\]
gives the cotrace $(a+d)(b+c)=\tau\tau'$, giving the type and genus of $H_2(\M)$.
Similarly, its Petrie length is determined by the trace
\[(a^2+bc)^2+(b(a+d))^2+(c(a+d))^2+(d^2+bc)^2\]
of $[x^2,y]$. Applying $H_{-1}$ or $D$ is achieved by replacing $x$ with
\[x^{-1}=\begin{pmatrix}d & b \\ c & a\end{pmatrix} \quad\hbox{or}\quad
z=\begin{pmatrix}c & a \\ d & b\end{pmatrix}.\]
Alternatively, one can replace $x$ with
\[z^{-1}=xy=\begin{pmatrix}b & a \\ d & c\end{pmatrix}\]
instead of $z$ when applying $D$.

The Frobenius triple-counting formula shows that $\Fr$ is the only map of type $\{3,7\}$ in $\O(G)$. For $x$ to correspond to such a map we require 
$ad+bc=1$, $a+d=t+1$, $t^2+1$ or $t^2+t+1$ and $b+c=1$. Without loss of generality we can use the solution
\[a=t, \quad b=1,\quad c=0,\quad d=t^2+1,\]
so that $\mathcal F$ is represented by the matrix
\[x=\begin{pmatrix}t & 1 \\ 0 & t^2+1\end{pmatrix}.\]
with trace $\tau=t^2+t+1$ and cotrace $\tau'=1$, confirming that the type is $\{3,7\}$ and the genus is $7$. Moreover, since $\sigma=t^2+1^2+0^2+(t^2+1)^2=t$, 
the Petrie length is $18$, so the extended type is $\{3,7\}_{18}$. This is therefore the map R7.1a in~\cite{Con09}. The dual map $D({\Fr})=$ R7.1b of extended type 
$\{7,3\}_{18}$ is represented by the matrix
\[z=\begin{pmatrix}0 & t \\ t^2+1 & 1\end{pmatrix}.\]

The map $H_2(\Fr)$ is represented by the matrix
\[\begin{pmatrix}t^2 & t^2+t+1 \\ 0 & t^2+t+1\end{pmatrix}\]
with trace $\tau^2=t+1$, cotrace $\tau\tau'=t^2+t+1$ and sum of squares $t^2+t$, so it has extended type $\{7,7\}_{18}$ and genus $55$.
It must therefore be one of the dual pair R55.52 in~\cite{Con09}; we will denote it by R55.32a, and its dual map $DH_2(\Fr)$, corresponding to the matrix
\[\begin{pmatrix}0 & t^2 \\ t^2+t+1 & t^2+t+1\end{pmatrix}\]
and also of type $\{7,7\}_{18}$, by R55.32b. (See Proposition~\ref{prop:SL_2(8)} for a proof that there are just two maps of type $\{7,7\}$ in $\O(G)$; 
the fact that the matrix $x^2$ for $H_2(\Fr)$ has distinct trace and cotrace shows that they are not self-dual.)

Iterating this process, we find that $H_2^2({\mathcal F})=H_4({\mathcal F})=H_3({\mathcal F})$ is represented by the matrix
\[\begin{pmatrix}t^2+t & t \\ 0 & t+1\end{pmatrix}\]
with trace $t^2+1$, cotrace $t$ and sum of squares $t^2+1$, so it has extended type $\{9,7\}_{14}$ and genus $63$.
It is therefore either R63.6b or R63.7b, since R63.6 and R63.7 are the only entries of this genus and extended type in~\cite{Con09}.
One can verify that it is R63.7b by checking that the matrices $x$ (given above) and $z$ corresponding to this map satisfy the defining relations for the full 
automorphism group of R63.7 given in~\cite{Con09}, but not those for R63.6. Specifically, the defining relations for R63.7, with generators $R$ and $S$ 
corresponding to our $x$ and $z$, and a third orientation-reversing generator $T$ inverting them both, are given as
\[R^{-7}=S^{-9}=(RS)^2=(S^{-1}R)^3=(RS^{-3}R^2)^2=1,\]
\[T^2=(RT)^2=(ST)^2=1.\]
One can check that the matrices
\[x=\begin{pmatrix}t^2+t & t \\ 0 & t+1\end{pmatrix} \quad\hbox{and}\quad z=\begin{pmatrix}0 & t^2+t \\ t+1 & t\end{pmatrix} \]
corresponding to $H_3(\Fr)$ satisfy these relations when substituted for $R$ and $S$. For example,
\[z^{-1}x=\begin{pmatrix}t & t^2+t \\ t+1 & 0\end{pmatrix}\begin{pmatrix}t^2+t & t \\ 0 & t+1\end{pmatrix}=
\begin{pmatrix}t^2+t+1 & t^2+1 \\ 1 & t^2+t\end{pmatrix}\]
has trace $1$, so $(z^{-1}x)^3=1$. On the other hand, the relations for R65.6 include $(RS^{-2}R)^2=1$, and since
\[xz^{-2}x=\begin{pmatrix}t^2+1 & t \\ t^2+t & t^2+t+1\end{pmatrix}\]
has trace $t\ne 0$ we have $(xz^{-2}x)^2\ne 1$. Thus $H_3({\mathcal F})$ is R63.7b.

\begin{proposition}\label{prop:SL_2(8)}
There are two orientably regular maps of type $\{7,7\}$ with orientation-preserving automorphism group $G=\SL_2(8)$. They form a dual pair.
\end{proposition}

\noindent{\sl Proof.} We will use the Frobenius triple-counting formula~(\ref{eq:Frobenius}), as in the proof of Proposition~\ref{prop:PSL_2(7)}, 
to count triples of type $(7,2,7)$ in $G$. There are three conjugacy classes $\mathcal X$, $\mathcal Z$ of elements of order $7$ in $G$, each 
containing  $2^3\cdot 3^2$ elements, and one class $\mathcal Y$ of $3^2\cdot 7$ involutions. Using the character values in~\cite{ATLAS}, 
we find that for each of the six choices of $\mathcal X\ne\mathcal Z$ the number of triples $(x,y,z)$ of type $(7,2,7)$ in $G$, with $x\in\mathcal X$ 
and $z\in\mathcal Z$, is
\[\frac{(2^3\cdot 3^2)\cdot(3^2\cdot 7)\cdot (2^3\cdot 3^2)}{2^3\cdot 3^2\cdot 7}
\left(1+\frac{1}{9}\sum_{j=0}^2(\zeta^{2^j}+\zeta^{-2^j})(\zeta^{2^{j+1}}+\zeta^{-2^{j+1}}) \right)\]
\[=\frac{(2^3\cdot 3^2)\cdot(3^2\cdot 7)\cdot (2^3\cdot 3^2)}{2^3\cdot 3^2\cdot 7}\cdot\frac{7}{9}=2^3\cdot 3^2\cdot 7\]
 where $\zeta=\exp(2\pi i/7)$, so the total number of such triples is
 \[2^4\cdot 3^3\cdot 7=2|{\rm Aut}\,G|.\]
These triples all generate $G$, since no maximal subgroup of $G$ contains  such a triple with non-conjugate $x$ and $z$, so there are two corresponding maps. 
 
 There are also three choices of classes $\mathcal X=\mathcal Z$, and the total number of corresponding triples in $G$ is
 \[3\cdot\frac{(2^3\cdot 3^2)\cdot(3^2\cdot 7)\cdot (2^3\cdot 3^2)}{2^3\cdot 3^2\cdot 7}
\left(1+\frac{1}{9}\sum_{j=0}^2(\zeta^{2^j}+\zeta^{-2^j})^2 \right)\]
\[=3\cdot\frac{(2^3\cdot 3^2)\cdot(3^2\cdot 7)\cdot (2^3\cdot 3^2)}{2^3\cdot 3^2\cdot 7}\cdot\frac{14}{9}
=2^4\cdot 3^3\cdot 7\]
Now $G$ has nine Sylow $2$-subgroups $T$, each with normaliser $N_G(T)\cong{\rm AGL}_1(8)\cong {\rm V}_8\rtimes {\rm C}_7$ generated by $48\cdot 7$ 
triples of type $(7,2,7)$, all with conjugate $x$ and $z$ (see Section~\ref{AGL_1(q)} for details). Since
\[9\cdot 48\cdot 7=2^4\cdot 3^3\cdot 7\]
these account for all the triples of type $(7,2,7)$ in $G$ with $x$ and $z$ conjugate; thus no such triples generate $G$, so they do not correspond to maps in $\O(G)$.

We have shown that there are two maps of type $\{7,7\}$ in $\O(G)$. One map corresponds to the pair $(t+1, t^2+1)$ of traces for $x$ and $z$, together with its 
Galois conjugates $(t^2+1, t^2+t+1)$ and $(t^2+t+1, t+1)$, while the other map corresponds to the pairs $(t+1, t^2+t+1)$, $(t^2+1, t+1)$ and $(t^2+t+1, t^2+1)$. 
Since the first three pairs differ from the last three by transposition of $x$ and $z$, these two maps form a dual pair. \hfill$\square$



\subsection{Summary of results for $\O(G)$}

Continuing in this way, using a similar criterion to identify the other maps of extended type $\{9,7\}_{14}$ or $\{7,9\}_{14}$, one eventually finds that there are fourteen maps 
$\M\in\O(G)$; they are described in Table~\ref{OSL28}, where the first four columns are analogues of those in Table~\ref{OPSL27} for $\PSL_2(7)$, with the Fricke--Macbeath 
map $\Fr$ replacing Klein's map $\K$. The final column gives the entries in~\cite{Con09} corresponding to the non-orientable quotients $\M/{\rm C}_2$ by the centre ${\rm C}_2$ 
of the full automorphism group $G\times{\rm C}_2$ of $\M$ (see Section~\ref{sec:non-orquot}).
 
 \begin{table}[ht]
\centering
\begin{tabular}{| p{2cm} | p{1.3cm} | p{2.8cm} | p{1.8cm} | p{1.2cm} |}
\hline
Entry in~\cite{Con09} & type & relationship to $\Fr$ & $H_2(\M)$ & $\M/{\rm C}_2$ \\
\hline\hline
R7.1a & $\{3,7\}_{18}$ & $\Fr$ & R55.32a & N8.1a\\
\hline
R7.1b & $\{7,3\}_{18}$ & $D(\Fr)$ & & N8.1b\\
\hline
R15.1a & $\{3,9\}_{14}$ & $H_2DH_3(\Fr)$ & R71.15a & N16.1a\\
\hline
R15.1b & $\{9,3\}_{14}$ & $DH_2DH_3(\Fr)$ &  & N16.1b\\
\hline
R55.32a & $\{7,7\}_{18}$ & $H_2(\Fr)$ & R63.7b & N56.5a\\
\hline
R55.32b & $\{7,7\}_{18}$ & $DH_2(\Fr)$ & R63.5b & N56.5b\\
\hline
R63.5a & $\{7,9\}_6$ & $DH_2DH_2(\Fr)$ & R71.15b & N64.3a \\
\hline
R63.5b & $\{9,7\}_6$ & $H_2DH_2(\Fr)$ & R63.6b & N64.3b\\
\hline
R63.6a & $\{7,9\}_{14}$ & $DH_3DH_2(\Fr)$ & R63.5a & N64.4a\\
\hline
R63.6b & $\{9,7\}_{14}$ & $H_3DH_2(\Fr)$ & R55.32b & N64.4b\\
\hline
R63.7a & $\{7,9\}_{14}$ & $DH_3(\Fr)$ & R15.1a & N64.5a\\
\hline
R63.7b & $\{9,7\}_{14}$ & $H_3(\Fr)=H_2^2(\Fr)$ & R7.1a & N64.5b\\
\hline
R71.15a & $\{9,9\}_{18}$ & $H_3DH_3(\Fr)$ & R63.7a & N72.9a \\
\hline
R71.15b & $\{9,9\}_{18}$ & $DH_3DH_3(\Fr)$ & R63.6a & N72.9b \\
\hline

\end{tabular}
\caption{The maps in $\O({\rm SL}_2(8))$}
\label{OSL28}
\end{table}

\begin{figure}[h!]
\begin{center}
\begin{tikzpicture}[scale=0.2, inner sep=0.8mm]

\node (a) at (-25,10) [shape=circle, fill=black] {};
\node (b) at (-15,10) [shape=circle, fill=black] {};
\node (c) at (-5,10) [shape=circle, fill=black] {};
\node (d) at (5,10) [shape=circle, fill=black] {};
\node (e) at (15,10) [shape=circle, fill=black] {};
\node (f) at (25,10) [shape=circle, fill=black] {};
\node (g) at (-10,4) [shape=circle, fill=black] {};
\node (h) at (10,4) [shape=circle, fill=black] {};
\node (i) at (-10,-4) [shape=circle, fill=black] {};
\node (j) at (10,-4) [shape=circle, fill=black] {};
\node (k) at (-15,-10) [shape=circle, fill=black] {};
\node (l) at (-5,-10) [shape=circle, fill=black] {};
\node (m) at (5,-10) [shape=circle, fill=black] {};
\node (n) at (15,-10) [shape=circle, fill=black] {};

\draw [thick, dashed] (a) to (b);
\draw [thick] (b) to (c) to (g) to (b);
\draw [thick, dashed] (c) to (d);
\draw [thick] (d) to (e) to (h) to (d); -
\draw [thick, dashed] (e) to (f);
\draw [thick, dashed] (g) to (i);
\draw [thick, dashed] (h) to (j);
\draw [thick] (i) to (k) to (l) to (i);
\draw [thick] (j) to (m) to (n) to (j);
\draw [thick, dashed, rounded corners] (k) to (-5,-12) to (m);
\draw [thick, dashed, rounded corners] (n) to (5,-12) to (l);

\node at (-25,12) {R7.1b};
\node at (-15,12) {R7.1a};
\node at (-5,12) {R63.7b};
\node at (5,12) {R63.7a};
\node at (15,12) {R15.1a};
\node at (25,12) {R15.1b};
\node at (-14.7,4.2) {R55.32a};
\node at (-14.7,-3.8) {R55.32b};
\node at (14.9,4.2) {R71.15a};
\node at (14.7,-3.8) {R71.15b};
\node at (-19.7,-9.8) {R63.5b};
\node at (19.7,-9.8) {R63.6a};
\node at (-3.5,-7.5) {R63.6b};
\node at (3.5,-7.5) {R63.5a};

\draw [thick]  (-13,10.5) to (-14.7,10) to (-13,9.5);
\draw [thick]  (13,10.5) to (14.7,10) to (13,9.5);
\draw [thick]  (-7,-10.5) to (-5.3,-10) to (-7,-9.5);
\draw [thick]  (7,-10.5) to (5.3,-10) to (7,-9.5);
\draw [thick]  (7,8.9) to (5.3,9.6) to (5.7,8);
\draw [thick]  (-7,8.9) to (-5.3,9.6) to (-5.7,8);
\draw [thick]  (13,-8.9) to (14.7,-9.6) to (14.3,-8);
\draw [thick]  (-13,-8.9) to (-14.7,-9.6) to (-14.3,-8);
\draw [thick]  (12,5.1) to (10.3,4.3) to (10.7,6.1);
\draw [thick]  (-12,5.1) to (-10.3,4.3) to (-10.7,6.1);
\draw [thick]  (8,-5.1) to (9.7,-4.3) to (9.3,-6.1);
\draw [thick]  (-8,-5.1) to (-9.7,-4.3) to (-9.3,-6.1);

\end{tikzpicture}

\end{center}
\caption{The graph $\O({\rm SL}_2(8))$, with names of maps} 
\label{fig:SL_2(8)}
\end{figure}
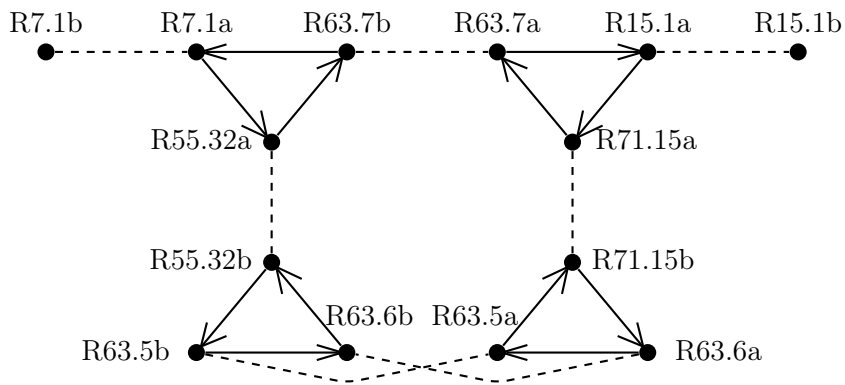

\begin{figure}[h!]
\begin{center}
\begin{tikzpicture}[scale=0.2, inner sep=0.8mm]

\node (a) at (-25,10) [shape=circle, fill=black] {};
\node (b) at (-15,10) [shape=circle, fill=black] {};
\node (c) at (-5,10) [shape=circle, fill=black] {};
\node (d) at (5,10) [shape=circle, fill=black] {};
\node (e) at (15,10) [shape=circle, fill=black] {};
\node (f) at (25,10) [shape=circle, fill=black] {};
\node (g) at (-10,4) [shape=circle, fill=black] {};
\node (h) at (10,4) [shape=circle, fill=black] {};
\node (i) at (-10,-4) [shape=circle, fill=black] {};
\node (j) at (10,-4) [shape=circle, fill=black] {};
\node (k) at (-15,-10) [shape=circle, fill=black] {};
\node (l) at (-5,-10) [shape=circle, fill=black] {};
\node (m) at (5,-10) [shape=circle, fill=black] {};
\node (n) at (15,-10) [shape=circle, fill=black] {};

\draw [thick, dashed] (a) to (b);
\draw [thick] (b) to (c) to (g) to (b);
\draw [thick, dashed] (c) to (d);
\draw [thick] (d) to (e) to (h) to (d); -
\draw [thick, dashed] (e) to (f);
\draw [thick, dashed] (g) to (i);
\draw [thick, dashed] (h) to (j);
\draw [thick] (i) to (k) to (l) to (i);
\draw [thick] (j) to (m) to (n) to (j);
\draw [thick, dashed, rounded corners] (k) to (-5,-12) to (m);
\draw [thick, dashed, rounded corners] (n) to (5,-12) to (l);

\node at (-25,12) {$\{7,3\}_{18}$};
\node at (-15,12) {$\{3,7\}_{18}$};
\node at ( -5,12) {$\{9,7\}_{14}$};
\node at (  5,12) {$\{7,9\}_{14}$};
\node at (15,12) {$\{3,9\}_{14}$};
\node at (25,12) {$\{9,3\}_{14}$};
\node at (-14.7,4.2) {$\{7,7\}_{18}$};
\node at (-14.7,-3.8) {$\{7,7\}_{18}$};
\node at (14.9,4.2) {$\{9,9\}_{18}$};
\node at (14.7,-3.8) {$\{9,9\}_{18}$};
\node at (-19,-9.8) {$\{9,7\}_{6}$};
\node at (19,-9.8) {$\{7,9\}_{14}$};
\node at (-3.5,-7.5) {$\{9,7\}_{14}$};
\node at (3.5,-7.5) {$\{7,9\}_{6}$};

\draw [thick]  (-13,10.5) to (-14.7,10) to (-13,9.5);
\draw [thick]  (13,10.5) to (14.7,10) to (13,9.5);
\draw [thick]  (-7,-10.5) to (-5.3,-10) to (-7,-9.5);
\draw [thick]  (7,-10.5) to (5.3,-10) to (7,-9.5);
\draw [thick]  (7,8.9) to (5.3,9.6) to (5.7,8);
\draw [thick]  (-7,8.9) to (-5.3,9.6) to (-5.7,8);
\draw [thick]  (13,-8.9) to (14.7,-9.6) to (14.3,-8);
\draw [thick]  (-13,-8.9) to (-14.7,-9.6) to (-14.3,-8);
\draw [thick]  (12,5.1) to (10.3,4.3) to (10.7,6.1);
\draw [thick]  (-12,5.1) to (-10.3,4.3) to (-10.7,6.1);
\draw [thick]  (8,-5.1) to (9.7,-4.3) to (9.3,-6.1);
\draw [thick]  (-8,-5.1) to (-9.7,-4.3) to (-9.3,-6.1);

\end{tikzpicture}

\end{center}
\caption{The graph $\O({\rm SL}_2(8))$, with types of maps} 
\label{fig:SL_2(8)types}
\end{figure}
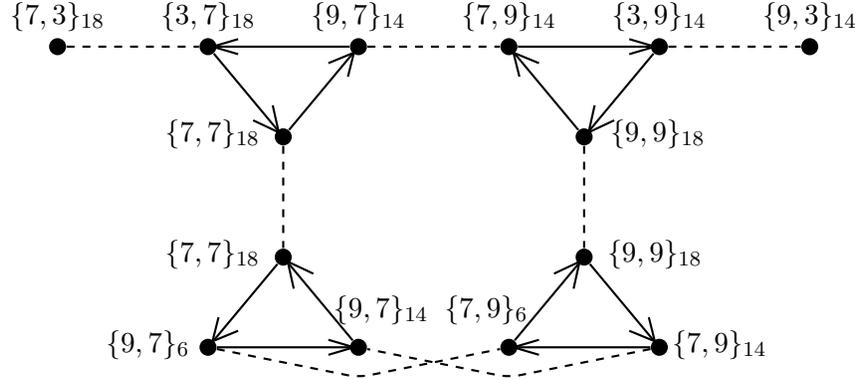

Figure~\ref{fig:SL_2(8)} shows the graph $\O(G)$, with broken and unbroken edges representing the actions of $D$ and $H_2$. The bilateral symmetry reveals an 
interesting duality between the parameters~$7$ and~$9$ appearing in the types of the maps in $\O(G)$ (see Figure~\ref{fig:SL_2(8)types}). However, this does 
not extend consistently to the Petrie lengths: for example, R7.1b of extended type $\{7,3\}_{18}$ is paired with R15.1b of extended type $\{9,3\}_{14}$, whereas R63.7b, 
of extended type $\{9,7\}_{14}$, is paired with R63.7a, of extended type $\{7,9\}_{14}$.

\medskip

\noindent{\bf Example} $H_3(\Fr)=H_2^2(\Fr)\cong{\rm R63.7b}$, of type $\{9,7\}$; this is confirmed by Figure~\ref{fig:H7}, which shows 
that the 3-holes of $\Fr$ have length $4+5=9$.


\subsection{Non-orientable quotients}\label{sec:non-orquot}

Each map $\M\in \O(G)$ is inner regular, with full automorphism group $G\times {\rm C}_2$: this is true for $\M=\Fr$ since the corresponding 
generators $x$ and $y$ are both inverted by the matrix 
\[\begin{pmatrix}t+1 & t \\ t & t+1\end{pmatrix}\in G;\]
now the graph $\O(G)$ is connected, and $\A\,\M$ is preserved by $D$ and $H_j$, so it is true for all $\M\in\O(G)$. The direct factor ${\rm C}_2$ 
reverses orientation, so each $\M$ has a non-orientable regular quotient $\overline\M=\M/{\rm C}_2$ with automorphism group $G$. Conversely,  
every map $\N\in\R(G)$ is non-orientable, since $G$, being simple, has no subgroup of index $2$, so $\N$ has a regular double cover $\M\in\O(G)$ 
with $\N=\overline\M$. This gives a bijection $\M\mapsto\N$ between the sets $\O(G)$ and $\R(G)$, where $\N$ has the same type as $\M$ (but 
half its Petrie length), and has genus $g+1$ if $\M$ has genus $g$. Since this bijection commutes with the operations $D$ and $H_j$, it induces an 
isomorphism $\O(G)\to\R(G)$ of directed graphs. Using their extended types to identify the maps $\N$ in the list of non-orientable regular maps 
in~\cite{Con09} gives the final column of Table~\ref{OSL28}.


\section{Some other Hurwitz groups} 

We have seen that if $G={\rm PSL}_2(7)$ or ${\rm SL}_2(8)$ then $\O(G)$ is connected, 
consisting of maps which are respectively all outer or all inner regular, with automorphism groups isomorphic to 
${\rm PGL}_2(7)$ or ${\rm SL}_2(8)\times{\rm C}_2$. 
In fact, a group $G$ can have both inner and outer regular maps in $\O(G)$, necessarily in different components. For instance, 
if $G={\rm PSL}_2(13)$ then $\O(G)$ contains three regular Hurwitz maps R14.1, R14.2 and R14.3, of genus $14$ and extended 
types $\{3,7\}_{12}$, $\{3,7\}_{26}$ and $\{3,7\}_{14}$. Calculations similar to those given earlier show that R14.2 is inner regular, 
with automorphism group $G\times{\rm C}_2$ and non-orientable regular quotient N15.1 of type $\{3,7\}_{13}$, while the other two 
maps are outer regular, with automorphism group ${\rm PGL}_2(13)$. This is interesting because, as shown by Streit~\cite{Str}, for 
each prime $p\equiv\pm 1$ mod~$(7)$ the three Hurwitz dessins in $\O({\rm PSL}_2(p))$ are Galois conjugate, under the Galois group 
${\rm C}_3$ of the real cyclotomic field $\Q(\cos 2\pi i/7)$, so although the orientation-preserving automorphism group of a dessin is a 
Galois invariant, the full automorphism group is not.

More generally, Wendy Hall~\cite{Hall77} has shown that a Hurwitz map $\M$ with ${\rm Aut}^+\M\cong{\rm PSL}_2(q)$ is inner regular if and only if 
$3-\tau^3$ is a square in $\F_q$, where $\tau$ is the trace of the canonical generator $x$ of order $7$. For primes $q=p\equiv\pm 1$ mod~$(7)$, where
 there are three Hurwitz maps, her examples $p=167$, $13$, $43$ and $181$ show that none, one, two or all three of them can be inner regular.

In order to realise alternating groups ${\rm A}_n$ as Hurwitz groups (Theorem~\ref{th:Con}), Conder~\cite{Con80} constructed a sequence 
of planar maps $\M_n$ of type $\{7,3\}$ which have Hurwitz maps ${\mathcal H}_n$ with ${\rm Aut}^+{\mathcal H}_n\cong{\rm A}_n$ as 
orientably regular covers. (Each $\M_n$ is presented as a coset diagram for a subgroup ${\rm A}_{n-1}<{\rm A}_n=\langle x, y\rangle$, but 
by shrinking the triangles for $x$ to points it can be interpreted as a cubic map with monodromy group ${\rm A}_n$.) The maps ${\mathcal H}_n$ 
are all outer regular: indeed, the bilateral symmetry in the construction of the maps $\M_n$ was designed to show that each map ${\mathcal H}_n$ 
has full automorphism group ${\rm S}_n$.

\medskip

\begin{figure}[htbp]
\begin{center}
\includegraphics[scale=0.2]{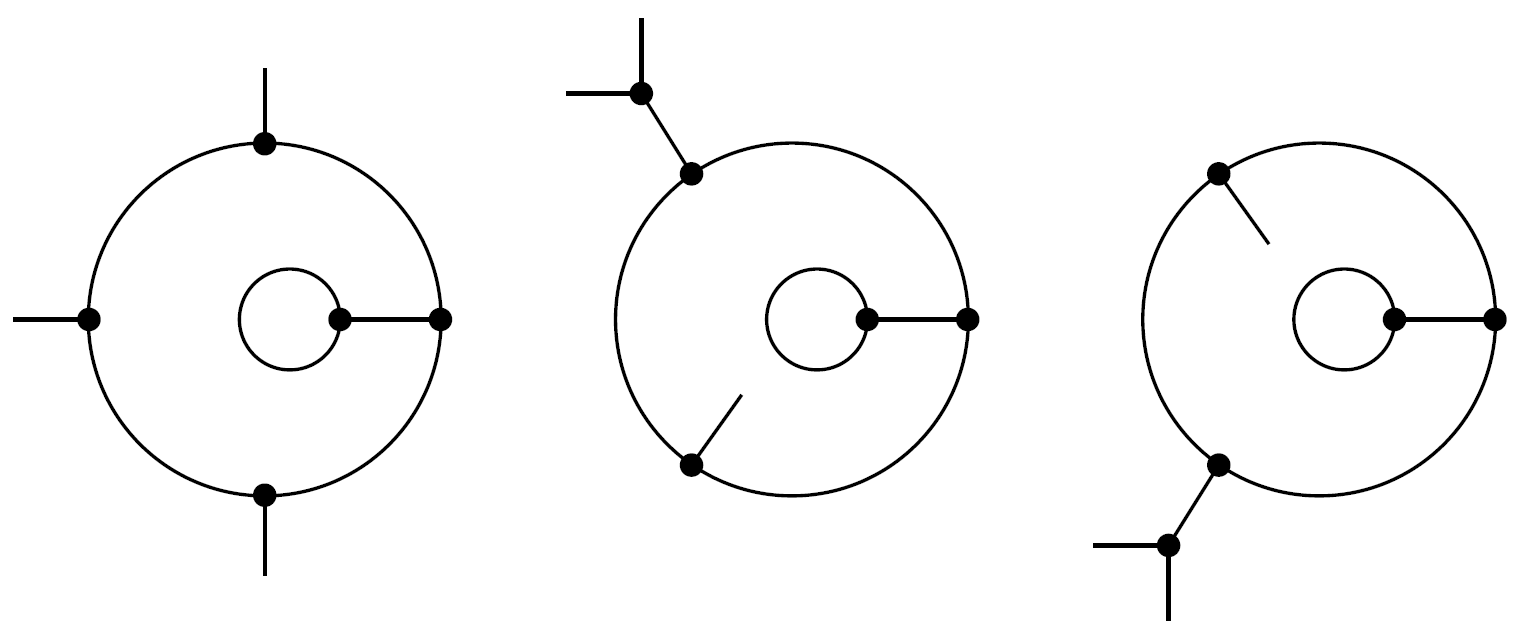}
\caption{\small Three maps of type $\{7,3\}$ with monodromy group ${\rm A}_{15}$.}
\label{fig:a_15}
\end{center}
\end{figure}

\noindent{\bf Example} Among the alternating groups, the smallest Hurwitz group is ${\rm A}_{15}$. There are three Hurwitz maps $\mathcal H$ 
with ${\rm Aut}^+{\mathcal H}\cong {\rm A}_{15}$. Instead of drawing them (they have genus $7\,783\,776\,001$), we show their planar quotients 
$\M={\mathcal H}/{\rm A}_{14}$ in Figure~\ref{fig:a_15}. On the left is the map $\M_{15}$, corresponding to Conder's diagram $B$ from his set 
of basic coset diagrams $A,\ldots,N$ in~\cite{Con80}; it is covered by the outer regular Hurwitz map ${\mathcal H}_{15}$. The other two maps are 
covered by a chiral pair of Hurwitz maps $\mathcal H$.

\medskip

By contrast, the simple Ree groups ${\rm Re}(q)={}^2G_2(q)$ ($q=3^e$ for odd $e\ge 3$) are also Hurwitz groups, but as shown in~\cite{Jon94} 
the Hurwitz maps associated with them are all chiral (the generator $x$ of order $3$ is not inverted by any automorphism). In the next section we 
will consider a more straightforward example of this last phenomenon, but this time not involving Hurwitz maps.
 
 
\section{${\rm AGL}_1(q), q=2^e$}\label{AGL_1(q)}
 
In this section, instead of considering individual groups such as $\PSL_2(7)$ and $\SL_2(8)$, we will consider an infinite family of groups $G$ 
for which $\O(G)$ exhibits uniform behaviour.

Let $G$ be the $1$-dimensional affine group ${\rm AGL}_1(q)$ for $q=2^e$, consisting of the affine transformations
\[t\mapsto at+b,\quad (a, b\in{\mathbb F}_q,\, a\ne 0)\]
of the field ${\mathbb F}_q$. This is a semidirect product $T\rtimes S$ of an elementary abelian normal subgroup
 \[T\cong({\mathbb F}_q,+)\cong{\rm V}_q=({\rm C}_2)^e\]
consisting of the translations $t\mapsto t+b$, by a complement
 \[S\cong({\mathbb F}^*_q,\times)\cong{\rm C}_{q-1},\]
consisting of the transformations $t\mapsto at, a\ne0$. To avoid trivial cases, we will assume from now on that $e\ge2$.
 
 \begin{theorem}
{\rm(a)} If $e>2$ there are $\phi(q-1)/e$ maps in $\M\in\O(G)$, all of type $\{q-1,q-1\}_4$ and genus $(q-1)(q-4)/4$; they are all chiral, 
and satisfy $D(\M)\cong H_{-1}(\M)$ and $H_2(\M)\cong \M$. -
\vskip2pt
\noindent{\rm(b)} If $e=2$ there is a single map $\M\in\O(G)$, the tetrahedral map $\{3,3\}$ of genus $0$, which is outer regular with 
${\rm Aut}\,\M\cong {\rm A\Gamma L}_1(4)\cong{\rm S}_4$.
\vskip2pt
\noindent{\rm(c)} For each $e\ge2$ the graph $\O(G)$ is connected.
 \end{theorem}
 
\noindent{\sl Proof.} 
The $q-1$ involutions in $T$ are all conjugate in $G$, and the remaining non-identity elements, all of order dividing $q-1$, form $q-2$ 
conjugacy classes of size $q$ (the non-identity cosets of $T$ in $G$). These fuse into $(q-2)/e$ orbits of size $qe$ under the action of
 \[{\rm Aut}\,G={\rm A\Gamma L}_1(q)\cong G\rtimes{\rm Gal}\,{\mathbb F}_q\cong G\rtimes{\rm C}_e.\]

Any map $\M\in\O(G)$ corresponds to a generating triple $(x,y,z)$ for $G$, where $y$ has order $2$ and $xT$ generates $G/T$. 
There are  $\phi(q-1)q$ choices for $x$ and $q-1$ choices for $y$, giving $\phi(q-1)q(q-1)=\phi(q-1)|{\rm Aut}\,G|/e$ triples, so 
there are $\phi(q-1)/e$ maps $\M$ in $\O(G)$. Since $x$ and $z$ have order $q-1$ these maps have type $\{q-1,q-1\}$ and hence 
have genus $(q-1)(q-4)/4$. Since $[x,y]$ is a non-identity element of $T$ it has order $2$, so the Petrie length is $4$. In particular, 
if $e=2$ there is a single map $\M\in\O(G)$; this is the tetrahedral map $\{3,3\}$, which is outer regular with ${\rm Aut}\,\M\cong
{\rm A\Gamma L}_1(4)\cong{\rm S}_4$. If $e>2$ then no element of ${\rm Gal}\,{\mathbb F}_q$ inverts ${\mathbb F}^*_q$, so 
the maps $\M\in\O(G)$ are all chiral; since $z$ is conjugate to $x^{-1}$ they satisfy $D(\M)\cong H_{-1}(\M)$. Since $x\mapsto x^2$
 is an automorphism of ${\mathbb F}_q$ they are all invariant under $H_2$. Since any two generators of $S$ are powers of each other, 
any two maps $\M, \M'\in\O(G)$ satisfy $\M'=H_j(\M)$ for some $j$ coprime to $q-1$, so $\O(G)$ is connected. (It is, in fact, a quotient 
of a Cayley graph for the group of units $U_{q-1}$.)  \hfill$\square$
 
 \medskip
 
These maps $\M$ are instances of the orientably regular embeddings of complete graphs $K_q$ constructed by Biggs in~\cite{Big71}, 
where he showed that $K_q$ has such an embedding if and only if $q$ is a prime power. The construction for odd $q$ is similar; 
see~\cite{JJ} for the classification of such maps.


\subsection{Examples} 

The maps $\M$ arising for $e=3$ are the chiral pair of Edmonds maps of type $\{7,7\}_4$, corresponding to the entry C7.2 in~\cite{Con09}; 
these both lie on the Fricke--Macbeath surface $\mathcal S$ of genus $7$ realising the Hurwitz group ${\rm SL}_2(8)$. The affine group 
${\rm Aut}\,\M\cong {\rm AGL}_1(8)$ is a subgroup of index $9$ in ${\rm Aut}\,{\mathcal S}={\rm Aut}^+\Fr\cong {\rm SL}_2(8)$: 
it is the stabiliser of $\infty$ in the natural representation of ${\rm SL}_2(8)$ on the projective line ${\mathbb P}^1({\mathbb F}_8)$, 
and also the normaliser of a Sylow $2$-subgroup. This inclusion lifts to an index~$9$ inclusion $\Delta(7,2,7)<\Delta=\Delta(7,2,3)$ of 
triangle groups, namely item~(B) in Singerman's list of triangle group inclusions~\cite{Sin72}.

For $e=4$ we have the chiral pair C45.2 of type $\{15,15\}_2$, and for $e=5$ we have the three chiral pairs C217.45--47 of type $\{31,31\}_4$. 
The index~$9$ inclusions of automorphism groups and triangle groups mentioned above for $e=3$ do not generalise to higher powers of $2$.

\begin{figure}[h!]
\begin{center}
\begin{tikzpicture}[scale=0.25, inner sep=0.8mm]

\node (a) at (5,-8.75) [shape=circle, fill=black] {};
\node (b) at (10,0) [shape=circle, fill=black] {};
\node (c) at (5,8.75) [shape=circle, fill=black] {};
\node (d) at (-5,8.75) [shape=circle, fill=black] {};
\node (e) at (-10,0) [shape=circle, fill=black] {};
\node (f) at (-5,-8.75) [shape=circle, fill=black] {};

\draw [thick] (a) to (b) to (c) to (d) to (e) to (f)to (a);
\draw [thick, dashed] (5,-8.35) to (-5,9.15);
\draw [very thick, dotted] (5,-9.15) to (-5,8.35);
\draw [thick, dashed] (10,0.2) to (-10,0.2);
\draw [very thick, dotted] (10,-0.2) to (-10,-0.2);
\draw [thick, dashed] (-5,-8.35) to (5,9.15);
\draw [very thick, dotted] (-5,-9.15) to (5,8.35);

\draw [thick] (3.5,-8.05) to (4.8,-8.75) to (3.5,-9.45);
\draw [thick] (9.8,-2.25) to (9.8,-0.25) to (8.2,-1.35);
\draw [thick] (5.2,6.5) to (5.2,8.5) to (6.8,7.4);
\draw [thick] (-3.5,8.05) to (-4.8,8.75) to (-3.5,9.45);
\draw [thick] (-9.8,2.25) to (-9.8,0.25) to (-8.2,1.35);
\draw [thick] (-5.2,-6.5) to (-5.2,-8.5) to (-6.8,-7.4);

\node at (7.5,-9) {$\M_1$};
\node at (12.5,0) {$\M_3$};
\node at (7.5,9) {$\M_5$};
\node at (-8,9) {$\M_{-1}$};
\node at (-13,0) {$\M_{-3}$};
\node at (-8,-9) {$\M_{-5}$};

\end{tikzpicture}

\end{center}
\caption{The graph $\O({\rm AGL}_1(32))$}
\label{fig:AGL_1(32)}
\end{figure}
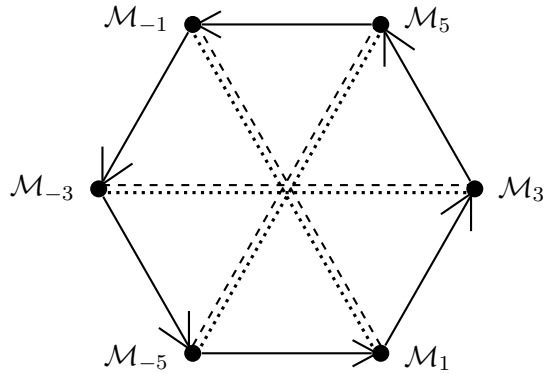

When $e=5$, if we take $a$ to be a generator of $\F_{32}^*=\langle a\mid a^{31}=1\rangle\cong{\rm C}_{31}$, then the chiral pairs of maps 
$\M_{\pm i}\;( i=1,3,5)$ correspond to the following mutually inverse pairs of orbits 
$\Omega_{\pm i}$ of ${\rm Gal}\,\F_{32}=\langle t\mapsto t^2\rangle\cong{\rm C}_5$ on $\F_{32}^*$:
\[\Omega_1=\{a^i\mid i=1, 2, 4, 8, 16\},\quad \Omega_{-1}=\{a^i\mid i=15, 30\equiv -1, 29, 27, 23\};\]
\[\Omega_3=\{a^i\mid i=3, 6, 12, 24, 17\},\quad \Omega_{-3}=\{a^i\mid i=7, 14, 28\equiv -3, 25, 19\};\]
\[\Omega_5=\{a^i\mid i=5,10, 20, 9, 18\},\quad \Omega_{-5}=\{a^i\mid i=11, 22, 13, 26\equiv -5, 21\}.\]
Then $D(\M_i)=H_{-1}(\M_i)=\M_{-i}$ and $H_2(\M_i)=\M_i$ for all $i$, while $H_3$ induces a $6$-cycle $(1,3,5,-1,-3,-5)$ on the subscripts $i$. 
The graph $\O({\rm AGL}_1(32))$ is shown in Figure~\ref{fig:AGL_1(32)}, with directed edges representing the action of $H_3$, and undirected 
dashed and dotted edges the actions of $D$ and $H_{-1}$. The identification of these maps with the entries C217.45--47 in~\cite{Con09} depends 
on the choice of $a$, or more precisely that of its minimal polynomial in the action on the additive group of $\F_{32}$, one of the six irreducible 
factors of the cyclotomic polynomial $\Phi_{31}(t)=t^{30}+t^{29}+\cdots+t+1$ in $\F_2[t]$ (see~\cite{JJ}). 

\section{Groups for $G$ which $\O(G)$ has many components}

In contrast with the groups $G={\rm AGL}_1(2^e)$, for which ${\mathcal O}(G)$ is connected for all $e$, we will now consider some families of groups $G$ 
for which ${\mathcal O}(G)$ has an unbounded number of connected components. The map operations $D$ and $H_j$ on orientably regular maps preserve 
the involution $y\in G$ in generating triples, up to the action of ${\rm Aut}\,G$, so if ${\rm Aut}\,G$ has $i=i(G)$ orbits on useful involutions $y\in G$, where 
`useful' means `member of a generating pair for $G$', then the number $c=c(G)$ of connected components of $\O(G)$ satisfies $c\ge i$.  The example $G={\rm S}_5$, 
with $i=2$ but $c=3$ (see Section~\ref{subsec:S5}), shows that $c$ can exceed $i$; the following example shows that both can be arbitrarily large.

\begin{lemma}\label{le:Sn}
If $G={\rm S}_n$ with $n\ge 5$ then each involution $y\in G$ is useful.
\end{lemma}

\noindent{\sl Proof.} It is sufficient to show that at least one involution in each conjugacy class of $G$ is useful. The result is straightforward if $y$ is a transposition 
$(i,j)$, since one can then take $x$ to be an $n$-cycle with $ix=j$, so we may assume that $y$ consists of $t$ transpositions and $n-2t$ fixed points for some $t\ge 2$.

First we will deal with the case $n\ge 8$. Let $m=\lfloor n/2\rfloor$, so that $m>3$. By the Bertrand--Chebyshev Theorem (see~\cite{SW}, for example) there is a prime 
$p$ such that $m<p<2m-2$, so $n/2<p\le n-3$. Let us take $x$ to have cycles of length $p$ and $n-p$ if $n$ is odd, and $p$, $n-p-1$ and $1$ if $n$ is even, so that $x$ 
is odd in either case. Now $x$ has at most three cycles, so given any $t\ge 2$ we can choose the involution $y$, which moves $2t\ge 4$ points, so that $H:=\langle x,y\rangle$ 
is transitive. It follows that $H$ must be primitive, for otherwise $H$ is contained in a wreath product ${\rm S}_a\wr{\rm S}_b$ for some proper factorisation $n=ab$ of $n$, 
which is impossible since $x$ has order divisible by $p$ whereas ${\rm S}_a\wr{\rm S}_b$ has order $(a!)^bb!$ coprime to $p$. Thus $H$ is a primitive group containing a 
$p$-cycle (a suitable power of $x$) with $n-p\ge 3$ fixed points, so by a classic theorem of Jordan (see~\cite[Theorem~13.9]{Wie64}) $H\ge{\rm A}_n$. Since $H$ contains 
the odd permutation $x$ we must have $H={\rm S}_n$, as required.

The case $n=5$ was considered in Section~\ref{subsec:S5}, while the cases $n=6$ and $n=7$ can easily be dealt with by hand or by using GAP.
\hfill$\square$

\medskip

The restriction $n\ge 5$ is required in this lemma: see Section~\ref{sec:small}.

\begin{corollary}
For each $m\in\mathbb N$ there exists $N_m\in\mathbb N$ such that if $n\ge N_m$ then $\O({\rm S}_n)$ has at least $m$ connected components.
\end{corollary}

\noindent{\sl Proof\/} By Lemma~\ref{le:Sn} the involutions in $G={\rm S}_n$ are all useful for $n\ge 5$. They have $\lfloor n/2\rfloor$ possible 
cycle-structures, so if $n\ne 6$ then ${\rm Aut}\,G\;(=G)$ has $\lfloor n/2\rfloor$ orbits on them. Thus $\O(G)$ has $c(G)\ge i(G)=\lfloor n/2\rfloor$ 
connected components for all $n\ge 7$. We may therefore take $N_m=\max\{2m, 7\}$. \hfill$\square$

\medskip

By adapting Conder's proof in~\cite{Con80} of Theorem~\ref{th:Con} one can also prove:

\begin{theorem}\label{th:An}
For each $m\in\mathbb N$ there exists $N'_m\in\mathbb N$ such that if $n\ge N'_m$ then $\O({\rm A}_n)$ has at least $m$ connected components containing Hurwitz maps.
\end{theorem}

\noindent{\sl Outline of proof\/} In~\cite{Con80} Conder proved that ${\rm A}_n$ is a Hurwitz group for each $n\ge 168$ (and also for some smaller $n$) by constructing 
coset diagrams for subgroups of index $n$ in $\Delta$, and showing that the induced permutation group on the cosets is ${\rm A}_n$. These diagrams are constructed by 
joining copies of $14$ basic cosets diagrams $A, B,\ldots, N$ of degrees between $14$ and $108$. An $i$-{\em handle} in a coset diagram for $\Delta$ is a pair $(a,b)$ of 
fixed points of $y$ with $a=bz^i$, where $i=1, 2$ or $3$. (In this paper we have transposed Conder's notation in~\cite{Con80} for the generators $x$ and $y$ of $\Delta$.) 
Two coset diagrams of degrees $d$ and $d'$ with $i$-handles $(a,b)$ and $(a',b')$ can be joined by an $i$-{\em join\/} to create a diagram of degree $d+d'$ by replacing 
$a$ and $a'$ with a transposition $(a,a')$ for $y$, and similarly for $b$ and $b'$.

Conder's costruction (simplified a little here) is as follows. His diagram $G$ of degree $42$ has three $1$-handles. First use $(1)$-joins to form a chain of $k$ copies of 
$G$ for some $k\ge 1$. For each of the $42$ congruence classes $[c]=[0], \ldots [41]$ mod~$(42)$, join two specified combinations of basic diagrams, of total degree 
$d_c\equiv c$ mod~$(42)$, to the ends of the chain, to give a coset diagram of degree $n=42k+d_c\equiv c$ mod~$(42)$. By taking $k=1, 2, \ldots$ this realises 
${\rm A}_n$ for each sufficiently large $n\in[c]$ as a quotient of $\Delta$ and hence as a Hurwitz group.

Each of the $k$ copies of $G$ in the chain has an unused $(1)$-handle, giving two fixed points of $y$. By using these to make $i$ further joins, for any $i\le\lfloor k/2\rfloor$, 
one can reduce the total number of fixed points of $y$ by $4i$; this does not change the fact that the resulting permutation group is ${\rm A}_n$ since Conder's proof of this 
still applies. We thus obtain $\lfloor k/2\rfloor+1$ Hurwitz maps in $\O({\rm A}_n)$, with mutually distinct cycle-structures for $y$ and hence in distinct components of this graph. 
Taking $k$ sufficiently large proves the result. \hfill$\square$


\section{The order of $\O(G)$}\label{sec:order}

One obvious problem which we have not yet addressed is to determine $|\O(G)|$ for any finite group $G$; however, a method for solving this is 
already known. The maps in $\O(G)$ correspond bijectively to the orbits of ${\rm Aut}\,G$ on generating pairs $x, y$ for $G$ satisfying $y^2=1$. 
Since ${\rm Aut}\,G$ acts semiregularly (i.e.~fixed-point-freely) on all generating sets for $G$, we have
\[|\O(G)|=\phi(G)/|{\rm Aut}\,G|,\]
where $\phi(G)$ is the number of such generating pairs $x, y$ for $G$. If $\sigma(G)$ denotes the total number of pairs $x, y\in G$ with $y^2=1$, 
then by applying Philip Hall's technique of M\"obius inversion in groups~\cite{Hall36} to the obvious equation
$\sigma(G)=\sum_{H\le G}\phi(H)$
we obtain
\begin{equation}\label{eq:phi}
\phi(G)=\sum_{H\le G}\mu(H)\sigma(H)
\end{equation}
where $\mu$ is the M\"obius function for the subgroup lattice of $G$, defined by
\[\mu(G)=1\quad\hbox{and}\quad \sum_{K\ge H}\mu(K)=0 \quad\hbox{for all}\quad H<G.\]
This applies to generating sets satisfying any given set of relations, but in our case, with $y^2=1$ the only relation, we have
\[\sigma(H)=|H|(|H|_2+1)\]
where $|H|_2$ is the number of involutions in $H$.

In~\cite{Hall36} Hall determined the function $\mu$ for many groups, including ${\rm PSL}_2(p)$ for primes $p$. Downs~\cite{Dow91} extended 
this to ${\rm PSL}_2(q)$ for prime powers $q=p^e$. Their results for $p>2$ are complicated but $|\O({\rm PSL}_2(13))|=33$ is a simple example. 
If $q=2^e$ then $|{\rm Aut}\,G|=e|G|$ and as in~\cite{DJ} we have
\[|\O(G)|=\frac{1}{e}\sum_{f|e}\mu\left(\frac{e}{f}\right)(2^f-1)(2^f-2)\]
where $\mu$ now denotes the M\"obius function of elementary number theory (see Sections~\ref{sec:small} and~\ref{sec:PSL_2(8)} for the cases 
$e=2$ and $e=3$). For any $q$ the sum in~(\ref{eq:phi}) is dominated by the summand with $H=G$, so that for $G={\rm PSL}_2(q)$ we have
\[|\O(G)|\sim\frac{|G|(|G|_2+1)}{|{\rm Aut}\,G|}\sim \frac{q^2}{e}\quad\hbox{or}\quad \frac{q^2}{4e}\quad\hbox{as}\quad q\to\infty,\]
where $q$ is respectively even or odd.

We close with a question: if $G={\rm PSL}_2(q)$, how does the number $c(G)$ of connected components of $\O(G)$ behave as $q\to\infty$?


\section*{Acknowledgement}
G\'{a}bor G\'{e}vay's research was supported by the Hungarian National Research, 
Development and Innovation Office, OTKA grant No.\ SNN 132625.


 \bigskip
 
 \noindent
Bolyai Institute\\
University of Szeged\\
H-6720 Szeged\\
Aradi v\'ertan\'uk tere 1\\
Hungary\\
{\tt gevay@math.u-szeged.hu}\\

\medskip

\noindent School of Mathematics\\
University of Southampton\\
Southampton SO17  1BJ\\
UK\\
{\tt G.A.Jones@maths.soton.ac.uk}\\

\end{document}